%% file: note.tex
\documentclass[a4paper,11pt]{article}%,dvipdfm

\usepackage{mathrsfs} 

\usepackage{amsthm}
\usepackage{hyperref}
\usepackage{amssymb}
\usepackage{latexsym}
\usepackage{amsmath}
\usepackage{amsfonts}
\usepackage{mathabx}
\usepackage[dvips]{graphicx}
\usepackage[utf8]{inputenc}

\usepackage[LGR,T1]{fontenc}%
\usepackage[french,english]{babel}
\usepackage{url}
\usepackage{enumerate}
\usepackage{hyperref}
\usepackage{xcolor}
\usepackage{fancyhdr}
\newtheorem{Theorem}{Theorem}[part]

\newtheorem{Proposition}{Proposition}[part]

\newtheorem{Lemma}{Lemma}[part]

\newtheorem{Remark}{Remark}[part]

\makeatletter \@addtoreset{equation}{section}

\@addtoreset{Definition}{section}

\@addtoreset{Theorem}{section}

\@addtoreset{Proposition}{section}

\@addtoreset{Property}{section}

\@addtoreset{Assumption}{section}

\@addtoreset{Corollary}{section}

\@addtoreset{Lemma}{section}

\@addtoreset{Remark}{section}

\@addtoreset{Example}{section}

\input{basecom}

\renewcommand{\whX}{X}

\definecolor{jcg}{HTML}{000000}

\definecolor{jfc}{RGB}{0,0,0} 

\definecolor{myred}{RGB}{0,0,0}

\title{A note on the $\cW_2$-convergence rate of 
the empirical measure of an ergodic $\R^d$-valued diffusion\footnote{This work is a companion paper to \cite{chassagneux2024computing}: A preliminary version of this note can be found in its first arXiv version, Section 5.} 
}

\author{
{\sc  
Jean-Francois Chassagneux \thanks{ENSAE-CREST and Institut Polytechnique de Paris, France. Email: {\tt chassagneux@ensae.fr}}
~\thanks{This work has been supported by the Chair ``Capital Markets Tomorrow: Modeling and Computational Issues'', a joint initiative of LPSM/Université Paris-Cité and Crédit Agricole CIB under the aegis of Institut Europlace de Finance.} 
\quad
Gilles Pag\`es} \thanks{Laboratoire de Probabilit\'es, Statistique et Mod\'elisation, UMR~8001, Sorbonne Universit\'e, case 158, 4, pl. Jussieu, F-75252 Paris Cedex 5, France. E-mail: {\tt  gilles.pages@sorbonne-universite.fr}}~\thanks{This research
benefited from the support of the ``Chaire Risques Financiers'', Fondation du Risque.}  
}
\begin{document}
\maketitle

\begin{abstract} 
    % \begin{center}
    %     \color{red}
    %     \LARGE Version 1.24
    % \end{center}
    % \bigskip
    In this note, we consider a Stochastic Differential Equation under a strong confluence and Lipschitz continuity assumption of the coefficients. For the unique stationary solution, we study the rate of convergence of its empirical measure toward the invariant probability measure. We provide rates for the Wasserstein distance, both in the mean quadratic and almost sure sense.
\end{abstract}

%\tableofcontents

\input{intro_note.tex}

\input{note_convtomu}

\small
\bibliographystyle{plain}
\bibliography{note_biblio}

\appendix
\section{Appendix}

\input{appendix}

\end{document}

%% file: basecom.tex
\newcommand{\cA}{\mathcal{A}}
\newcommand{\cB}{\mathcal{B}}

\newcommand{\cF}{\mathcal{F}}

\newcommand{\cL}{\mathcal{L}}
\newcommand{\cM}{\mathcal{M}}

\newcommand{\cO}{\mathcal{O}}
\newcommand{\cP}{\mathcal{P}}

\newcommand{\cW}{\mathcal{W}}

\newcommand{\E}{\mathbb{E}}
\newcommand{\F}{\mathbb{F}}

\renewcommand{\P}{\mathbb{P}}

\newcommand{\R}{\mathbb{R}}

\newcommand{\dX}
    {\ensuremath{\delta X}}

\def \proof{{\noindent \bf Proof. }}
\def \eproof{\hbox{ }\hfill$\Box$}

\newcommand{\ud}{\mathrm{d}}

\newcommand{\HYP}[1]
    {\ensuremath{({H#1} ) }}

\newcommand{\1}{{\bf 1}}
    
\newcommand{\set}[1]
    {\ensuremath{\{ #1 \}}}
    
\newcommand{\HP}[1] %L2DPDT sur 0,T
    {\ensuremath{\mathscr{H}^{#1}}}

\newcommand{\esp}[1]{\ensuremath{\mathbb{E} \! \left[#1\right] }}

\renewcommand{\Xi}[1]{X_{i #1}}

\newcommand{\whX}{\widehat{X}}
\newcommand{\hm}{\hat{\mu}}
\newcommand{\hn}{\hat{\nu}}

\newcommand{\Lip}[1]{\ensuremath{[#1]_{\rm Lip}}}

\def \a{\alpha}

\def \d{\delta}
\def \la{\lambda}
\def \ve{\varepsilon}

\def \s{\sigma}

\newcommand{\ent}[1]{\ensuremath{ \lfloor #1 \rfloor}}

\def \hypl2{\HYP{^\star}}

%% file: intro_note.tex
%\section{Rate of convergence of the empirical measure toward the stationary measure}
\section{Introduction}
\label{se rate conv emp meas to stationary distrib}

In this note, we study the convergence of the empirical measure of a  Stochastic Differential Equation (SDE)  toward its stationary distribution. This work is a companion paper to \cite{chassagneux2024computing}.  This question has received some interest in the past years. Several recent works ,see e.g. \cite{wang2024sharp,trevisan2024wasserstein,wang2024asymptotics}, study the convergence in $L^q$ norm for the $p$-Wasserstein distance in some compact abstract spaces. The article also \cite{lei2020convergence} presents  some concentration results, see also the references therein. 
In \cite{chassagneux2024computing}, our motivation is to obtain rates of convergence for a numerical method computing the stationary distribution of a McKean-Vlasov SDE by ergodic simulation. Related to this question, recent papers \cite{du2023empirical,du2023self}  are studying the convergence of the empirical measure of self-interacting diffusion (SID) toward the stationary distribution of some related McKean-Vlasov SDEs. A key point in proving the ergodic behavior of the SID and its Euler scheme, and in establishing a convergence rate, is first to control the rate of convergence in Wasserstein distance between the empirical measure of a standard SDE and its invariant distribution. 

To be more precise, 
let $(\Omega,\cA,\P,\F:=(\cF_t)_{t\ge 0})$ be  a filtered probability space, with $\F$ satisfying the usual conditions. Let $q,d$ be two positive integers. We are  given a $q$-dimensional $\F$-Brownian motion $W$ (in particular independent of $\cF_0$) and 
we consider 
the following SDE
\begin{align}\label{eq solution sde hat}
        \ud \whX_t= B(\whX_t)\ud t + \Sigma(\whX_t) \ud W_t,
    \end{align}
where the function $(B,\Sigma):\R^d \rightarrow \R^d\times\cM_{d,q}(\R)$ satisfies the following  Lipschitz continuity assumption
   \begin{align}
       \HYP{L}\!:& \quad| B(x) - B(y) |  + \| \Sigma(x) - \Sigma(y)\|_{_F}
       \le L  |x-y|   ,
   \end{align}
   for a positive constant $L$.
For later use, we will denote by $\cL$ the infinitesimal generator of the above diffusion.\\
We will also assume  that the coefficients $(B,\Sigma)$ satisfy the following strong confluence assumption: there exists some $\alpha > 0$
\begin{align}\label{eq strong confluence for classic SDE}
    \HYP{C}\!:& \quad 2(B(x)-B(y)|x-y) + \| \Sigma(x)-\Sigma(y)\|^2_{_F} \le -\alpha |x-y|^2, \quad \text{ for $\alpha>0$.}
\end{align}

\begin{Lemma}\label{le sigma bounded exp moments}
Under \HYP{C} and \HYP{L}, there exists a unique stationary distribution  $\hat{\mu}$ to \eqref{eq solution sde hat} and 
\[
        \exists\, a>0\; \mbox{ such that }\; \int|y|^{2+a}  \hm(dy)<+\infty.
        \]
If, moreover, $\Sigma$ is bounded then 
\[
    \exists\, \lambda>0\; \mbox{ such that }\; \int e^{\lambda |y|^{2}}  \hm(dy)<+\infty.
\]     
\end{Lemma}
\noindent For the sake of completeness, a proof of this Lemma is given in the appendix.

\vspace{2mm} \noindent
 From now on, we assume that  $\whX_0\sim \hat \mu$ and thus $\whX_t \sim \hat{\mu}$ for all $t \ge 0.$
\medskip
Let us denote, for $t\ge 0$,
\begin{align}\label{eq empirical meas classical sde}
    \hn_t(dy) = \frac 1t \int_0^t \delta_{\whX_s}ds
\end{align}
the empirical measure associated to the stationary solution $\whX$ of~\eqref{eq solution sde hat}.  

Our main results concern the convergence of the  $(\hat\nu_t)_{t\ge 0}$ 
% given by
% \begin{align}
%     \nu_t := \frac1t\int_0^t \delta_{X_s} \ud s, t > 0 \text{ and $\nu_0 = \delta_{X_0}$, }
% \end{align}
towards the stationary distribution $\hat{\mu}$. They are summed up  in the two following Theorems.
% \textcolor{red}{!! Ces deux théorèmes doivent être un invariant des différentes versions de cette note car ils sont cités dans le papier numérique !!}
\begin{Theorem}\label{th conv L2 as minimal ass}
    Assume \HYP{L} and \HYP{C} are in force. \textcolor{jfc}{Let $a>0$ be given by Lemma \ref{le sigma bounded exp moments}.}
        % Under~\eqref{eq strong confluence for classic SDE}, we know that
        % \[
        % \exists\, a>0\; \mbox{ such that }\; \int|y|^{2+a}  \hm(dy)<+\infty.
        % \]
        %
        \begin{enumerate}[i)]
    \item Then, we have that 
        \begin{align}\label{eq pr conv L2}
            \big(\E {\cal W}^2_2(\hn_t, \hm)\big)^\frac12 = O\left( t^{-\frac{a}{2(2(d+2)+a(d+3))}} \right).
        \end{align}
    \item If moreover $\Sigma$ is uniformly elliptic \footnote{A matrix valued function $\R^d \ni x \mapsto \sigma(x) \in \cM_{d,q}$ is \emph{uniformly elliptic} if it satisfies, for some $\varsigma_0 > 0$:
    \begin{align}\label{eq de unif elliptic}
        \upsilon^\top \sigma \sigma^\top(x) \upsilon \ge \varsigma_0|\upsilon|^2, \quad \forall x \in \R^d, \forall \upsilon \in \R^d.  
    \end{align} }, 
    % and  that
    %     \begin{align}\label{eq moment ass as}
    %         \exists\, a >0\; \mbox{ such that }\; \int|y|^{2+a}  \hm(dy)<+\infty.    
    %     \end{align}
        \color{jfc}
        then it holds that
        \begin{align}\label{eq pr conv L2 upgraded}
            \big(\E {\cal W}^2_2(\hn_t, \hm)\big)^\frac12 = O\left( t^{-\frac{a}{2(d+2)(a+2)}} \right),
        \end{align}
        and $\P-a.s.$
        \begin{align}\label{eq pr conv as}
            \cW_2(\hn_t,\hm) = o_\eta\big(t^{-\zeta_a}\log^{\frac12+\eta}(t)\big)\,, \quad \mbox{ for every $\eta > 0$},   
        \end{align}
        where \textcolor{jfc}{$ \zeta_a =\frac{a^2}{2(2d+(a+2)(d+2)(a+d+2))}$}. 
        \color{black}
    \end{enumerate}
    \end{Theorem}

    \vspace{2mm}
    Under a boundedness assumption on the diffusion coefficient $\Sigma$, leading to bounded exponential moment for the invariant distribution, we can improve the previous rate of convergence.

    \begin{Theorem} \label{th conv L2 as exp moment}
        Let \HYP{L} and \HYP{C} hold and
              assume that $\Sigma$ is bounded.
              % \textcolor{red}{devrait suffire pour l'utilisation qu'on en fait dans le papier numérique}
        % \begin{align} \label{eq ass exp moment L2}
        %     \exists\, \lambda>0\; \mbox{ such that }\; \int e^{\lambda |y|^2}  \hm(dy)<+\infty ,  
        % \end{align}
        % e.g. because $\Sigma$ is bounded (see Lemma~\ref{le sigma bounded exp moments}).
        \begin{enumerate}[i)]
        \item Then, it holds that 
        \begin{align}\label{eq pr conv L2 exp moment}
            \big(\E\, \cW_{2}^2(\hn_t,\hm)\big)^\frac12 = O \left(t^{-\frac1{2(d+3)}} 
            (\log t)^{\frac{d+2}{2(d+3)}}\right).
        \end{align}
        \item If moreover $\Sigma$ is uniformly elliptic\footnotemark[1],
        \color{jfc}
        then it holds that
        \begin{align}\label{eq pr conv L2 exp moment upgraded}
            \big(\E\, \cW_{2}^2(\hn_t,\hm)\big)^\frac12 = O \left(t^{-\frac1{2(d+2)}} 
            (\log t)^{\frac12}\right),
        \end{align}
and $\P-a.s.$ 
        \begin{align}\label{eq pr conv as exp moment}
            \cW_{2}(\hn_t,\hm)  = t^{-\frac1{2(d+2)}}o_\eta\left(\log(t)^{\frac12+\eta}\right), \quad \mbox{ for every $\eta > 0$}.
            %O \left(
            %(\log t)^{\frac{d+2}{d+3}}t^{-\frac1{2(d+3)}} \right).
        \end{align}
        \color{black}
    \end{enumerate}
\end{Theorem}

\begin{Remark}
    The rates in Theorem \ref{th conv L2 as exp moment} are obtained, up to a log term, by sending $a$ towards infinity in the rate given in Theorem \ref{th conv L2 as minimal ass}.
\end{Remark}

The main part of this note is devoted to the proof of the above Theorems. 
Theorem \ref{th conv L2 as minimal ass}(i) is proved in Section \ref{se L2 conv}, see Proposition \ref{pr conv L2}(a). Theorem \ref{th conv L2 as minimal ass} (ii)  is proved in Section \ref{se as conv}, see Proposition \ref{pr W2 as conv}.
Theorem \ref{th conv L2 as exp moment}(i) is proved in Section \ref{se L2 conv}, see Proposition \ref{pr conv L2}(b).  Theorem \ref{th conv L2 as exp moment}(ii) is proved in Section \ref{se as conv}, see Proposition \ref{pr W2 as conv exp moment}. %\textcolor{red}{ pour le cas L2 on n'a pas besoin de $\Sigma$ bornée mais seulement de moment exp, verifier ce qui se passe pour a.s. sans soute dans la preuve concentraiton pour martingale}

\medskip 
\noindent 
% Our approach to obtain the previously convergence rates relies on the existence of {\em coboundaries}   of  Lipschitz functions whose existence and useful properties are  established in the next lemmas. 

The starting point of our proof, see Section \ref{se L2 conv}, relies on  the classical smoothing procedure also adopted in \cite{du2023empirical} or \cite{fathi2013transport}. But, in a second phase, we rely on a different approach based on coboundary functions (solution to the Poisson equation). In particular, this new point of view allows us to obtain almost sure rates of convergence for the Wasserstein distance which are not present in \cite{du2023empirical}.

% \textcolor{red}{
% It is important to notice that to prove the ergodic behavior of the SID, one must first control the rate of convergence in Wasserstein distance between the empirical measure of a standard SDE and its invariant measure. For this specific result, we retrieve and improve, under some further reasonnable assumptions, the rate obtained 
% in \cite{du2023empirical} in our setting (natural empirical measure). As a first step, our proof, see Section \ref{se rate conv emp meas to stationary distrib}, follows  the classical smoothening procedure also adopted in \cite{du2023empirical} or \cite{fathi2013transport}. But, in a second phase, we rely on a different approach based on coboundary functions (solution to the Poisson equation). In particular, this new point of view allows us to obtain almost sure rates of convergence for the Wasserstein distance which are not present in \cite{du2023empirical}.
% }
% \textcolor{red}{Section \ref{se rate conv emp meas to stationary distrib} is concerned with the rate of convergence of the empirical measure toward the stationary measure for classical SDEs: it is a key part to obtain the global error bounds once applied to the stationnary MKV SDEs.} 

The rest of the paper is organised as follows. Section \ref{se coboundaries and prelim} presents preliminary results on coboundary functions. Section \ref{se L2 conv} contains the proof for the $L^2$ convergence of the Wasserstein distance. Section \ref{se as conv} is concerned with the almost sure convergence. Finally, we gather in the appendix the proofs of some useful technical Lemmas.

\color{black}

\paragraph{Notations}

$\bullet$ $(u\,|\,v) = \sum_{1\le i \le d}u^i v^i$ denotes the canonical inner product on $\R^d$  of $u= (u^1,\ldots, u^d)$, $v=(v^1, \ldots,v^d)\!\in \R^d$, $|\cdot|$ the associated Euclidean norm. 

$\bullet$ The set $\cM_{d,q}(\R)$ is the set of real matrices with $d$ rows and $q$ columns. For a matrix $M$, its transpose matrix is denoted $M^\top$. 
The quantity $(A\,|\, B)_{_F}= \sum_{ij} a_{ij}b_{ij} = \mathrm{Tr}[A^\top B]$ denotes the ``Frobenius'' inner product between  matrices $A=[a_j]$ and $B=[b_{ij}] \! \in  \cM_{d,q}(\R)$ (set of matrices with $d$ rows and $q$ columns), $\| \cdot\|_{F}$ the associated Frobenius norm.

$\bullet$ \color{jfc} For a real bounded continuous function $f$, we denote by $\|f\|_{\sup}$ its supremum norm.\color{black}

$\bullet$ For a Lipschitz continuous function $f$, we denote by $\Lip{f}$ its Lipschitz constant w.r.t. the canonical Euclidean norm.

$\bullet$ A function $f:\R^d \rightarrow \R$ is in $C^{2,1}(\cO)$, for an open subset $\cO \subset \R^d$, if it is twice differentiable with second order derivatives which are Lipschitz continuous on $\cO$.

$\bullet$ For an $\R^d$-valued random variable $X$, we denote $\| X\|_p := (\esp{|X|^p})^\frac1p$, for $p>0$.

% A matrix valued function $\R^d \ni x \mapsto \sigma(x) \in \cM_{d,q}$ is \emph{uniformly elliptic} if it satisfies, for some $\epsilon > 0$:
% \begin{align}
%     \upsilon^\top \sigma \sigma^\top(x) \upsilon \ge \epsilon|\upsilon|^2, \quad \forall x \in \R^d, \forall \upsilon \in \R^d.  
% \end{align} 

%% file: note_convtomu.tex
\section{About coboundaries and preliminary rate results}

\label{se coboundaries and prelim}

As mentioned above, we rely, in our method of proof, on coboundary functions and some of their properties.  Our first result is to build coboundaries for Lipschitz functions when $\Sigma$ is uniformly elliptic in the sense of \eqref{eq de unif elliptic}.

\begin{Proposition} \label{le unif elliptic cobord}
 \noindent \color{jfc} Assume that \HYP{L} and \HYP{C} are in force. Let $f:\R^d \rightarrow \R$ be a  Lipschitz function. 
 The function $\omega_f$ defined 
    \begin{equation}\label{eq:cobord lip f}
        \omega_f(x) = \int_0^{+\infty}\big(P_tf(x)-\hm(f)\big)dt\;,\,\quad x \in \R^d,
        \end{equation}
    where $(P_t)_{t \ge 0}$ is the semi-group associated to \eqref{eq solution sde hat}, is Lipschitz-continuous with 
    \begin{align}\label{eq lip const wf first}
    \Lip{\omega_f} \le \frac{2 \Lip{f}}{\alpha}.    
    \end{align}
    \color{black}
    If, moreover, $\Sigma$ is uniformly elliptic, then $\omega_f$ is a $C^{2,1}(\R^d)$ solution of
    \begin{align}
        %(B(x)|\nabla \omega_f(x)) + \frac12(\partial^2_{xx} \omega_f(x)|\Sigma \Sigma^\top(x))_{_F} = \hm(f)-f(x) 
        %\;,\; x \in \R^d
        \cL\omega_f= \hm(f)-f 
        \;
    \end{align}
    % and, when $f$ is bounded, 
    % \begin{align}\label{eq lip const wf second}
    % \Lip{\omega_f} \le  c \|f\|_{\sup},
    % \end{align}
    % for some  constant $c>0$. 
    and  it also holds that 
    \begin{align}\label{eq expansion cobound along emp meas}
       \frac1t \int_0^t \set{f(\whX_s)-\hat\mu(f)} \ud s = \frac{\omega_f(\whX_0)-\omega_f(\whX_t)}{t}  + \frac1t \int_0^t \partial_x \omega_f(\whX_s)\Sigma(\whX_s)\ud W_s.
    \end{align}
\end{Proposition}

%{   
 %   \color{blue} 
    \proof 
    
   \textbf{Step 1.} Under the confluence condition~\eqref{eq strong confluence for classic SDE}, we know that  for every $x,y\!\in \R^d$ and every $t\ge 0$,
\begin{align}\label{eq basic inequality}
  \|X^x_t-X^y_t\|_2 \le |x-y|e^{-\frac{\a}2 t}.  
\end{align}
As a consequence for any $\mu \!\in {\cal P}_1(\R^d)$, with obvious notations
\[
\forall\, t\ge 0, \quad \|X^x_t-X_t^{\mu}\|_{1}   \le e^{-\frac{\a}2 t} \bigg( \int|x-y|\mu(dy)\bigg)
\]
owing to the Markov property at $t=0$. In particular, 
\begin{equation}\label{eq:Pgeom}
\forall\, t\ge 0, \quad |P_tf(x) -\mu(f)| \le [f]_{\rm Lip} \|X^x_t-X_t^{\mu}\|_{1} \le [f]_{\rm Lip}e^{-\frac{\a}2 t} \int|x-y|\mu(dy) ,
\end{equation}
so that $\omega_f$ is well defined by \eqref{eq:cobord lip f} and   Lipschitz continuous since 
    \begin{equation}\label{eq:Lipvarphi}
    |\omega_f(x)-\omega_f(x')| \le \int_0^\infty|P_tf(x)-P_tf(x')|dt \le   [f]_{\rm Lip}|x-y| 
    \int_0^\infty e^{-\frac{\a}2 s}\ud s =  \frac{ 2 [f]_{\rm Lip}}{\a}|x-y| .
    \end{equation}
 
\textbf{Step 2.} To deal with the differentiability of $\omega_f$, we call upon results from PDE theory. 
Now let $\cB:=\cB_1(x_0)$ be the open ball centered at an arbitrary $x_0 \in \R^d$. We define $w$ the solution to the following elliptic PDE 
 
\begin{align} \label{eq elliptic pde on B}
    \cL w= \hat{\mu}(f)-f \text{ on } \cB \text{ and } w=\omega_f \text{ on } \partial \cB.
\end{align}
From \cite[Theorem 6.13]{gilbarg1977elliptic}, the function $w \in C^{2,1}(\cB)\cap C^0(\bar{\cB})$. For a given $x \in \bar{\cB}$, we introduce the stopping time (w.r.t. the completed filtration of $W$) 
\begin{align*}
    \tau_x = \inf \set{t \ge 0| X^x_t \in \partial B }
\end{align*}
From \cite[Chapter II]{bass1998diffusions}, we know that $\P(\tau_x < +\infty)=1$ and that the following probabilistic representation holds true:
\begin{align}\label{eq rep w}
    w(x) = \esp{\omega_f(X^x_{\tau_x}) + \int_0^{\tau_x}\big(f(X^x_s)-\hat{\mu}(f)\big)\ud s }.
\end{align}
On the other hand, it follows from  the strong Markov property and the definition of $\omega_f$ that
\begin{align*}
    \omega_f(X^x_{\tau_x}) &= \int_0^{+\infty}\big(P_tf(X^x_{\tau_x})-\hm(f)\big)dt
    =\int_0^{+\infty}\big(\esp{f(X^x_{\tau_x+t})|\cF_{\tau_x}}-\hm(f)\big)dt
    \\
    &=
    \esp{\int_{\tau_x}^{+\infty} \big( f(X^{x}_t) -\hat{\mu}(f)\big) \ud s| \cF_{\tau_x}}.
\end{align*}
Inserting the previous equality into~\eqref{eq rep w}, we obtain
\begin{align}
    w(x) &= \esp{\int_0^{+\infty}\big(f(X^x_s)-\hat{\mu}(f)\big)\ud s }
    =\omega_f(x).
\end{align}
Thus  $\omega_f \in C^{2,1}(\cB) $ and satisfies~\eqref{eq elliptic pde on B}. Finally, observe that $x_0$ is arbitrary so that $\omega_f \in C^{2,1}(\R^d) $. 
%\tiny \textcolor{magenta}{One can get estimates for second order derivatives - boundedness in fact - with a dependence on the ellipticity constant obviously} \normalsize

\noindent Applying Ito's formula, we get 
\begin{align}
\nonumber     \omega(\whX_t) &= \omega(\whX_0) + \int_0^t \cL  \omega_f(\whX_s) \ud s + \int_0^t \partial_x \omega(\whX_s)\Sigma(\whX_s)\ud W_s \\
    &= \omega(\whX_0) + \int_0^t \set{\hat{\mu}(f)-f(\whX_s)} \ud s + \int_0^t \partial_x \omega(\whX_s)\Sigma(\whX_s)\ud W_s
\end{align}
which concludes the proof for this step. \eproof\\
    
%}
\color{black}

\vspace{2mm}
The above Proposition gives an expansion of the error $\hat{\nu}_t(f)-\hat{\mu}(f)$ in the almost sure sense. In the following Lemma, we obtain a control in expectation that allows us to remove the uniform ellipticity condition.
\begin{Lemma}\label{le only lip but no ps} 
    Assume that \HYP{L} and \HYP{C} are in force.
    There exists a  positive constant $C$, such that, for all $f:\R^d \rightarrow \R$  Lipschitz continuous,
    for all $t \ge 0$, 
    \begin{align}
        \esp{\Big|\frac1t\int_0^t (f(\whX_s)-\hat{\mu}(f)) \ud s\Big|} \le C\frac{\Lip{f}}{\sqrt{t}}. 
    \end{align}
\end{Lemma}

%{
%    \color{blue}
    \proof Let $\varepsilon \in (0,1)$ and $B$ a Brownian motion independent from $W$. Introduce the SDE
\begin{align}
    \ud X^\ve_t = B(X^\ve_t)\ud t + \Sigma(X^\ve_t) \ud W_t + \varepsilon \ud B_t.
\end{align}
Since condition~\eqref{eq strong confluence for classic SDE} is in force, we know from classical arguments that there is a unique stationary measure $\mu^\varepsilon \in \cP_2(\R^d)$ for this process. We also assume that $[X_0^\ve] = \mu^\ve$.
\\
We compute, denoting $\dX = X^\ve - \whX$,
\begin{align}
    e^{\a t}|\dX_t|^2 &= |\dX_0|^2 + \a \int_0^te^{\a s}|\dX_s|^2\ud s +\ve^2\int_0^te^{\a s}\ud s
    \\
    &+ 
        2\int_0^te^{\a s}\dX_s\left(\Sigma(X^\ve_s)-\Sigma(\whX_s)\right) \ud W_s
        + 2 \ve \int_0^te^{\a s}\dX_s \ud B_s
        \label{eq starting point aux results}
        \\
        &+\int_0^t e^{\a s}[2(B(X^\ve_s)-B(\whX_s)|\dX_s) + \| \Sigma(X^\ve_s)-\Sigma(\whX_s)\|^2_{_F}] \ud s. \nonumber
\end{align}
Using~\eqref{eq strong confluence for classic SDE}, we obtain 
\begin{align}
    |\dX_t|^2 &\le e^{-\a t}|\dX_0|^2 + \frac{\ve^2}{\a}+
        2\int_0^te^{\a (s-t)}\dX_s\left(\Sigma(X^\ve_s)-\Sigma(\whX_s)\right) \ud W_s
        + 2 \ve \int_0^te^{\a (s-t)}\dX_s \ud B_s.
\end{align}
Taking expectation in the above inequality we obtain, letting $t \rightarrow +\infty$, that 
\begin{align}
    \cW_{2}^2(\hm,\mu^\ve) \le \frac{\ve^2}{\a}.
\end{align}
Classical arguments also show that, for all $T \ge 0$, there exists $C_T$ such that
\begin{align}\label{eq same path}
    \esp{\sup_{t \in [0,T]} |\delta X_t|^2} \le C_T \ve^2.
\end{align}

\noindent  We observe that 
\begin{align*}
        \frac1t\int_0^t (f(\whX_s)-\hm(f)) \ud s = 
        \frac1t\int_0^t (f(\whX_s)-f(X^\ve_s)) \ud s
        +
        \frac1t\int_0^t (f(X^\ve_s)-\mu^\ve(f)) \ud s
        + \set{\mu^\ve-\hm}(f)
\end{align*}
Since the dynamics of $X^\varepsilon$ is uniformly elliptic, we can invoque 
Proposition~\ref{le unif elliptic cobord} to get 
\begin{align*}
    \Big|\frac1t\int_0^t (f(\whX_s)-\hm(f)) \ud s\Big|
    &\le 
    \frac{\Lip{f}}{t}\int_0^t |\dX_s|\ud s 
    + 2\frac{\Lip{f}}{\a}\frac {|X^\ve_t-X^\ve_0|}t + \Lip{f}\frac{\ve}{\sqrt{\a}} 
    \\
    &+\frac1t \Big|\int_0^t \partial_x \omega^\ve_{f}(X^\ve_s) \Sigma(X^\ve_s) \ud W_s \Big| + \ve \frac1t \Big|\int_0^t \partial_x \omega^\ve_{f}(X^\ve_s)  \ud B_s\Big|,
\end{align*}
where $\omega^\ve_f$ is the coboundary function introduced in Proposition~\ref{le unif elliptic cobord} applied with $X^\varepsilon$.
We compute that 
\begin{align}
    \esp{\Big|\int_0^t \partial_x \omega^\ve_{f}(X^\ve_s) \Sigma(X^\ve_s) \ud W_s \Big|^2}^\frac12 \le C\Lip{f} \sqrt{t},
\end{align}
since $\sup_{t\ge 0}\esp{|X^\ve_s|^2} = \int y^2 \mu^\ve(\ud y)\le C$. Importantly, $C$ can be chosen independent of $t \ge 0$ and $\ve \in (0,1)$.
Similarly, we obtain that 
\begin{align}
    \esp{|X^\epsilon_t-X^\epsilon_0|} \le C \sqrt{t}\; \text{ and } \;\esp{\Big|\int_0^t \partial_x \omega^\ve_{f}(X^\ve_s)  \ud B_s\Big|} \le  C\Lip{f} \sqrt{t}.
\end{align}
Thus, combining the previous inequality with~\eqref{eq same path}, we have 
\begin{align*}
    \esp{\Big|\frac1t\int_0^t (f(\whX_s)-\hm(f)) \ud s\Big|}
    \le 
C \Lip{f}\left(\ve + t^{-\frac12}\right),
\end{align*}
and the proof is concluded by letting $\ve \rightarrow 0$.
\eproof
%}

\vspace{2mm}
\color{jfc}
The following Lemma states an alternative control of the Lipschitz constant of the coboundary function when the source term  is bounded and the diffusion coefficient is uniformly elliptic.

\begin{Lemma}\label{le lip wf unif ellip bounded}
    Assume that \HYP{L} and \HYP{C} are in force and that $\Sigma$ is uniformly elliptic. Let $f:\R^d \rightarrow \R$ be a  Lipschitz continuous and bounded function.  Then $\omega_f$ defined in \eqref{eq:cobord lip f} satisfies 
    \begin{align}\label{eq lip const wf second}
    \Lip{\omega_f} \le  c \|f\|_{\sup},
    \end{align}
    for some  constant $c>0$. Moreover, it holds that, for all $t > 0$,
    \begin{align}
        \esp{\Big|\frac1t\int_0^t (f(\whX_s)-\hat{\mu}(f)) \ud s\Big|} \le C\frac{\|f\|_{\sup}}{\sqrt{t}}.
    \end{align}
\end{Lemma}
\proof \textbf{Step 1.} For $\epsilon>0$, let $B^\epsilon$ and $\Sigma^\epsilon$ be mollification of the coefficients $B$ and $\Sigma$. We observe that $\HYP{L}$, $\HYP{C}$ are still in force for these mollified coefficients: We denote by $\hat{\mu}^\epsilon$ the stationary measure. We also observe that $\Sigma^\epsilon$ is uniformly elliptic with ellipticity constant $\frac{\varsigma_0}{2}$ for $\epsilon$ small enough (have in mind \eqref{eq de unif elliptic}). Let $(P^\epsilon_t)$ the semi-group associated to the diffusion with mollified coefficients and define
\begin{equation}\label{eq:cobord mollified lip f}
        \omega^\epsilon_f(x) = \int_0^{+\infty}\big(P^\epsilon_tf(x)-\hm^\epsilon(f)\big)dt\;,\,\quad x \in \R^d.
        \end{equation}
We obtain, from classical computations, see e.g. the beginning of the proof of Lemma \ref{le only lip but no ps} , that 
\begin{align}
    \lim_{\epsilon \rightarrow 0} \omega^\epsilon_f(x) = \omega_f(x)
\end{align}
which implies that $
\Lip{\omega_f} \le \liminf_\epsilon\Lip{\omega^\epsilon_f}.$ The rest of this step is dedicated to proving that $\Lip{\omega^\epsilon_f} \le c \|f\|_\infty$ for $c>0$ that does not depend on $\epsilon$.  

\color{black}
\textbf{Step 2.} From now on we drop the reference to the mollified coefficient for notational simplicity and we assume directly that $B,\Sigma$ are smooth. Let \(X^x\) denote the solution of SDE starting from \(x\), and let
\(J_t^x\) be the derivative flow, namely
\[
    J_t^x v = \lim_{\eta \rightarrow 0} \frac{X^{x+\eta v}_t - X^{x}_t}{\eta},
    \qquad v\in\mathbb{R}^d .
\]  
From \eqref{eq basic inequality}, we get that
\begin{align}\label{eq decroissance J}
    \mathbb{E}\big[|J_t^x v|^2\big]
    \leq e^{-\alpha t}|v|^2 .    
\end{align}

 Put
\[
    a(x):=\Sigma(x)\Sigma(x)^\top 
\]
and since $\Sigma$ is uniformly elliptic, \(a(x)\geq\varsigma_0 I_d\), \(a(x)\) is invertible and
\begin{align}\label{eq unif ellip}
    \big|\Sigma(x)^\top a(x)^{-1} z\big|^2
    =
    z^\top a(x)^{-1}z
    \leq \varsigma_0^{-1}|z|^2, \quad z \in \R^d .
\end{align}

  Let \(h:[0,t]\to\mathbb{R}\) be absolutely continuous with \(h(0)=0\) and
\(h(t)=1\). The Bismut-Elworthy-Li formula gives, for every
\(v\in\mathbb{R}^d\),
\[
    \nabla P_t f(x)\cdot v
    =
    \mathbb{E}\left[
        f(X_t^x)
        \int_0^t
        \dot h_s
        \left(
            \Sigma(X_s^x)^\top a(X_s^x)^{-1}J_s^x v
            \,\middle|\, dW_s
        \right)
    \right].
\]
By Cauchy-Schwarz inequality, Itô isometry, \eqref{eq decroissance J}, and  \eqref{eq unif ellip}, we obtain
\[
\begin{aligned}
    |\nabla P_t f(x)\cdot v|
    &\leq
    \|f\|_\infty
    \left(
        \mathbb{E}
        \int_0^t
        |\dot h_s|^2
        \big|
            \Sigma(X_s^x)^\top a(X_s^x)^{-1}J_s^x v
        \big|^2
        ds
    \right)^{1/2}                                                        \\
    &\leq
    \frac{ \|f\|_\infty}{\sqrt{\varsigma_0}}
    \left(
        \int_0^t
        |\dot h_s|^2
        \mathbb{E}|J_s^x v|^2
        ds
    \right)^{1/2}                                                        \\
    &\leq
    \frac{ \|f\|_\infty |v|}{\sqrt{\varsigma_0}}
    \left(
        \int_0^t
        |\dot h_s|^2 e^{-\alpha s}
        ds
    \right)^{1/2}.
\end{aligned}
\]

We now choose \(h\) depending on \(t \in (0,\infty)\). First, for $t\le 1$, we take
\(\dot h_s=t^{-1}\) on \([0,t]\) and  we get
\begin{align}\label{eq control small time}
    |\nabla P_t f(x)|
    \leq
    \frac{ \|f\|_\infty}{\sqrt{\varsigma_0}}
    \left(
        \frac{1}{t^2}\int_0^t e^{-\alpha s}\,ds
    \right)^{1/2}
    \leq
    C\frac{\|f\|_\infty}{\sqrt t},
    \qquad 0<t\leq 1.
\end{align}
    
For $t > 1$, we choose instead
\[
    \dot h_s
    =
    \frac{2}{t}\mathbf{1}_{[t/2,t]}(s).
\]
Then \(h(0)=0\), \(h(t)=1\), and
 
\begin{align}
    |\nabla P_t f(x)|
    &\leq
    \frac{ \|f\|_\infty}{\sqrt{\varsigma_0}}
    \left(
        \frac{4}{t^2}\int_{t/2}^t e^{-\alpha s}\,ds
    \right)^{1/2}  \nonumber
    \\
    &\leq
    C \|f\|_\infty\, t^{-1} e^{-\alpha t/4}. \label{eq control large times}
\end{align}

Combining \eqref{eq control small time} and \eqref{eq control large times},   we obtain
the global estimate

\begin{align}\label{eq global estimate}
    \sup_{x\in\mathbb{R}^d}|\nabla P_t f(x)|
    \leq
    C\|f\|_\infty \left(\frac{1}{\sqrt t}1_{[0,1]}(t)+e^{-\alpha t/4}\right),
    \qquad t>0. 
\end{align}
Since the above bound is integrable on $(0,\infty)$, we obtain that
\[
    \nabla\omega_f(x)
    =
    \int_0^\infty \nabla P_t f(x)\,dt .
\]
Integrating \eqref{eq global estimate}, we do get 
\[
\begin{aligned}
    |\nabla\omega_f(x)|
    % &\leq
    % \int_0^\infty |\nabla P_t f(x)|\,dt                          \\
    % &\leq
    % C\|f\|_\infty
    % \int_0^\infty t^{-1/2}e^{-\theta t}\,dt                       \\
    &\leq
    c\|f\|_\infty , \text{ for some $c>0$},
\end{aligned}
\]
which concludes the proof for this step.
\\
\color{jfc}
\textbf{Step 3.} From \eqref{eq expansion cobound along emp meas}, we obtain 
\begin{align*}
  \esp{\left|\frac1t \int_0^t \set{f(\whX_s)-\hat\mu(f)} \ud s\right|} &\le \esp{ \left|\frac{\omega_f(\whX_0)-\omega_f(\whX_t)}{t} \right|} + \esp{ \left|\frac1t \int_0^t \partial_x \omega_f(\whX_s)\Sigma(\whX_s)\ud W_s\right|},
 \\
 &\le \frac{c\|f\|_{\sup}}{t}\esp{|X_t-X_0|} +
   \frac1t\esp{ \int_0^t \left| \partial_x \omega_f(\whX_s)\Sigma(\whX_s)\right|^2\ud s}^\frac12.
\end{align*}
We observe that $\esp{|X_t-X_0|}\le C\sqrt{t}$ and 
$ \esp{ \int_0^t \left| \partial_x \omega_f(\whX_s)\Sigma(\whX_s)\right|^2\ud s}^\frac12 \le C \|f\|_{\sup}\sqrt{t}$. This yields
\[
\esp{\left|\frac1t \int_0^t \set{f(\whX_s)-\hat\mu(f)} \ud s\right|} \le C\frac{\|f\|_{\sup}}{\sqrt{t}}
\]
and concludes the proof for this step.
\eproof
\color{black}

\vspace{2mm}
\noindent We conclude this section by stating two Lemmas, which are also useful in the proof of our main results.

% whose proofs are given in the Appendix~\ref{subsec expo moment} and Appendix~\ref{subsec hoeffding like} respectively, for sake of completeness.

\begin{Lemma}\label{le control SDE}
   Assume  that \HYP{L} and \HYP{C} 
   %the strong confluence condition~\eqref{eq strong confluence for classic SDE} 
   are in force. %\textcolor{red}{@Gil \& JFC: donner nom d'hypo}. 
   %Checker si on dit qu'alors moment d'ordre 2}.
    
    \smallskip
    \noindent $(a)$  Then,  
    \begin{align}
        \esp{\frac1t|\whX_t - \whX_0|} = O(t^{-1}) \quad\mbox{ as }t\to +\infty.
    \end{align} 
$(b)$ Moreover, assume that $\sup_{t\ge 0}\E\, |\whX_t|^{2+a} <+\infty$ for some $a\!\in [0, 2]$. Then for every $\eta >0$
    \begin{align}\label{eq: rateX_t/tas}
\frac1t |\whX_t - \whX_0| = o_\eta\big(t^{-\frac{1+a}{2+a}}(\log t)^{\frac 1{1+a/2}+\eta}\big)\quad\mbox{ as }\quad t\to +\infty.
    \end{align} 

%    Assume moreover, that $\int|y|^{4}  \hm(dy)$. Then, for $\epsilon > 0$,
%    \begin{align}  
%        \frac1t|\whX_t - \whX_0| = o\left(t^{-\frac34}\log(t)^{\frac14+\epsilon}\right) \; \text{ a.s.}
%    \end{align}
\end{Lemma}
\color{black}
%\proof 
\noindent {\em Proof}.  
$(a)$ The first statement follows the fact that, under~\eqref{eq strong confluence for classic SDE}, $ \E\, |\whX_t| \le \int|y|   \hm(dy)<+\infty$.

\noindent $(b)$ Let us consider the $C^2$-function $\Psi:\R^d\to \R$ defined by $\Psi_a(x)= (1+|x|^2)^{\frac12 + a/4}$. One checks that 
\[
\forall\, x\!\in \R^d,\qquad\nabla \Psi(x) = \big(1+\frac a2\big)\frac{x}{(1+|x|^2)^{\frac12-a/4}}\quad \mbox{and} \quad \nabla^2\Psi(x)= (1+\frac a2)\frac{(1+|x|^2)^{1/2}I_d-xx^{\top}}{(1+|x|^2)^{1-a/4}}.
\]
% It follows from~\eqref{eq strong confluence for classic SDE} that there exists  $\a' >  0$ and $K' \ge 0$ such that
% \[
% \forall\, x\!\in \R^d, \qquad  (x\,|\, B(x))+\tfrac 12 \|\s(x)\|_{F}^2 \le K'-\a' |x|^2
% \]
Since \eqref{eq strong confluence for classic SDE} is in force, \eqref{eq Hajek for classic SDE} holds too.
 Then, It\^o's formula applied to $e^{\la' \cdot}\Psi(\whX)$ with $\la' = (1+\frac a2)\a'$, yields
\begin{align*}
e^{\la' t}\Psi(\whX_t)& =\Psi(\whX_0)  +\int_0^te^{\la' s} \big(H_s\,|\,\ud W_s\big)
\\
&+\int_0^t  e^{\la' s}\Big[ \a'(1+\frac a2) \Psi(\whX_s)+ (1+\frac a2)\frac{(\whX_s\,|\,B(\whX_s))}{(1+|\whX_s)|^2)^{\frac12-a/4}}+\tfrac 12 {\rm Tr}\big(\Sigma\Sigma^{\top}(\whX_s)\nabla^2\Psi (\whX_s)\big)\Big]\ud s
%\\  &\qquad\qquad
\end{align*}
where $\displaystyle H_t=(1+\frac a2) \frac{\s(\whX_t)^{\top}\whX_t}{(1+|\whX_t|^2)^{\frac12-a/4}}$. Now 
\begin{align*}
{\rm Tr}\big(\Sigma\Sigma(\whX_t)^{\top}\nabla^2\Psi(\whX_t)\big) &
=(1+\frac a2)\bigg(\frac{\|\Sigma(\whX_s)\|^2_{F}}{(1+|\whX_s|^2)^{\frac12-a/4}}-\frac{|\Sigma(\whX_t)\whX_t|^2}{(1+|\whX_s|^2)^{1-a/4}}\bigg) \le (1+\frac a2)\frac{\|\Sigma(\whX_s)\|^2_{F}}{(1+|\whX_s|^2)^{\frac12-a/4}}
\end{align*}
so that,  
\begin{align*} 
    (1+\frac a2)\frac{(\whX_s\,|\,B(\whX_s))}{(1+|\whX_s|^2)^{1/2-a/4}}+\tfrac12 {\rm Tr}\big(\Sigma\Sigma(\whX_t)^{\top}\nabla^2\Psi(\whX_t)\big) 
&\le(1+\frac a2)\frac{(\whX_s\,|\,B(\whX_s))+\tfrac12 \|\Sigma(\whX_s)\|^2_{F}}{(1+|\whX_s|^2)^{\frac12-a/4}}\\
&\le(1+\frac a2) \frac{\beta'-\a' |\whX_s|^2}{(1+|\whX_s|^2)^{\frac12-a/4}}\\
&\le (1+\frac a2)( \a'+K' -\a' \Psi(\whX_s)).
%\le  \beta+\a'-\a'(1+ |\whX_t|^2)
\end{align*}
Hence
\[
\Psi(\whX_t) \le e^{-\la' t} \Psi(\whX_0) +\frac{\a'+K'}{\a'} + e^{-\la' t}\int_0^te^{\la' s} (H_s\,|\,\ud W_s).
\]
The integrand $H_t$ satisfies
\[
%\frac{\s(\whX_s)^{\top}\whX_s}{(1+|\whX_s|^2)^{1/2}}
|H_t| \le \frac{\|\Sigma(\whX_t)\|_{_F} |\whX_t|}{(1+|\whX_t|^2)^{\frac12-a/4}}\le \frac{C(1+|\whX_t|^2)}{(1+|\whX_t|^2)^{\frac12-a/4}}= C(1+|\whX_t|^2)^{\frac12+a/4},
\]
so that $\sup_{t\ge 0}\E\, |H_t|^2 \le C' (1+\sup_{t\ge 0}\E\, |\whX_t|^{2+a}) <+\infty$. 
\textcolor{black}{
Hence, mimicking the proof of claim~$(b)$ of Proposition 3.1 in \cite{chassagneux2024computing},}
 one derives using Stochastic Kronecker's Lemma, see Lemma \ref{le kronecker gil} in the Appendix, that, for every $\eta>0$, $\Psi(\whX_t)= o_{\eta}\big(t^{1/2}(\log t)^{1+\eta} \big)$ so that 
\[
\hskip 3,75cm \frac{|\whX_t|}{t}\le \frac{\Psi(\whX_t)^{\frac 2{2+a}}}{t}  =o_{\eta}\big(t^{-\frac{1+a}{2+a}}(\log t)^{\frac{1}{1+a/2}+\eta} \big).\hskip 3,75cm\Box
\]
% the stochastic integral is a true square integral martingale. 

\color{black}

This last result is concerned with concentration inequalities for Brownian martingales; its proof is given in the Appendix for sake of completeness. 

 \begin{Lemma}\label{le hoeffding like}
    Let $(Z_t)$ be a   progressively  measurable process and let  $\psi: (0,+\infty) \to (0,+\infty)$.\\
 $(a)$  If $(Z_t)$ is bounded then for every $\ell>0$, we have %\textcolor{green}{@JCF : \`a ce stade on pourrait ne garder que $\psi(t)$}
    \begin{align}\label{eq thanks hoeff}
         \P\left(\Big|\int_0^t Z_s \ud W_s\Big| > \ell \psi(t)\right) \le 
         2 e^{-\frac{\ell^2}{2|Z|_{\infty}^2}t^{- 1}\psi(t)^2}.
     \end{align}

\noindent $(b)$ If $\sup_{t\ge0}\esp{|Z_t|^{2+a}}\le C$, for some $a>0$. Then, the following holds 
    \begin{align} \label{eq encore merci bdg}
        \P\left(\Big| \int_0^tZ_s \ud W_s\Big| > \ell \psi(t)\right) \le 
        C \frac{t^{{\frac{a+2}2}}}{\ell^{a+2}\psi(t)^{a+2}}.
    \end{align}
 \end{Lemma}
 
%{
 %   \color{blue}
%}

\section{$L^2$ convergence}

\color{black}
\label{se L2 conv}
In this section, we obtain $L^2$-convergence rates for the Wasserstein distance under various integrability assumptions on the stationary distribution. In particular, the result for the case of exponential moments is new, to the best of our knowledge.  Though the approach by regularisation that we follow here may look rather classical, it has a some novelty as we rely on coboundary functions introduced in Section \ref{se coboundaries and prelim}.
The proofs given here should also be seen as a gentle introduction to those of the next section, which deals with the almost sure convergence rate.

\color{black}

\vspace{2mm}
\noindent Before starting let us recall two classical properties of the (quadratic) Wasserstein distance, for every $\mu, \, \tilde \mu ,\, m \!\in {\cal P}_2(\R^d)$, 
\begin{equation}\label{eq:ident1} \quad {\cal W}_2(\mu, \mu*m) \le \bigg( \int_{\R^d} |\xi|^2 m(\ud \xi)\bigg)^{1/2}
\end{equation}
and if $\mu= g.\lambda_d$ and $\tilde \mu= \tilde g.\lambda_d$  then it follows from \cite[Lemma 2.2]{horowitz1994mean},
%(see~\cite{}, \textcolor{red}{cf. papier cit\'e chez les chinois. En fait c'est juste une  majo de type TV!!}) 
\begin{equation}\label{eq:ident2 new}
 {\cal W}_2^2(\mu,\tilde \mu)\le 3 \int_{\R^d} |\xi|^2 |g(\xi)-\tilde g(\xi)| \ud \xi.
\end{equation}

\medskip

\color{black}
\noindent The starting point of the main proofs is a classical smoothening of the problem which is summarised in the following small technical Lemma.
\begin{Lemma} \label{le smoothing generic}
    Let $\rho $ be a positive mollifier with support   $[-1,1]^d$ and let $\zeta$ be a   $\rho(y)\lambda_d(dy)$-distributed random vector. Let $\ve>0$ and let  us define   
 $$
 \rho_{\ve}(u) = \ve^{-d}\rho\Big(\frac{u}{\ve}\Big).
 $$
Note that $\Lip{\rho_{\ve}} = \ve^{-(d+1)}\Lip{\rho}=  \ve^{-(d+1)}\|\nabla \rho\|_{\sup}$ and $\|\rho_{\ve} \|_{\sup} \le \ve^{-d}\|\rho \|_{\sup}$. 
\\
Setting 
\begin{equation}\label{eq:dfnphieps le}
\phi_{\ve}(\xi, x) =\rho_{\ve}(\xi-x)-\int_{\R^d}\rho_{\ve}(\xi-y)\hm(dy) 
\end{equation}
it holds that, for $t\ge 0$, 
\begin{align}\label{eq:decomp0 le}
    {\cal W}_2(\hn_t, \hm)^2
   &\le 6\ve^2 \|\zeta\|^2_2 + 9 \int_{\R^d} |\xi|^2  \Big|\frac1t\int_0^t \phi_{\ve} (\xi, \whX_s) ds \Big| \ud \xi.   
   \end{align}
\end{Lemma}
\color{black}
\proof It follows from the triangle inequality that
\begin{align*}
{\cal W}_2(\hn_t, \hm) &\le {\cal W}_2(\hn_t, \hn_t*[\ve \zeta]) +  {\cal W}_2(\hm, \hm *[\ve\zeta])+{\cal W}_2(\hn_t *[\ve\zeta], \hm *[\ve\zeta]) .
\end{align*}
Moreover, one checks that  both  distributions $\hn_t*[\ve \zeta]$ and $ \hm *[\ve\zeta]$ have densities respectively  given for every $t\ge 0$ by 
\[
\frac{d\hn_t *[\ve\zeta]}{d\lambda_d} =   \frac{1}{t}\int_0^t\rho_{\ve}(\xi-\whX_s) \ud s \quad\mbox{ and }\quad  
\frac{d\hm  *[\ve\zeta]}{d\lambda_d} = \int_{\R^d}\rho_{\ve}(\xi-y)\hm(dy).  
\]
Then \eqref{eq:ident2 new} yields that  
\begin{align*} %\label{eq:decomp0 as unbounded}
 {\cal W}_2(\hn_t, \hm)^2
&\le 6\ve^2 \|\zeta\|^2_2 + 9 \int_{\R^d} |\xi|^2  \Big|\frac1t\int_0^t \phi_{\ve} (\xi, \whX_s) ds \Big| \ud \xi,   
\end{align*}
which concludes the proof.
\eproof
\color{black}

%\color{yellow}
\begin{Proposition}\label{pr conv L2}
    Assume \HYP{C} and \HYP{L} hold.
    $(a)$ Let $a>0$ from Lemma \ref{le sigma bounded exp moments} such that
    \[
      \int|y|^{2+a}  \hm(dy)<+\infty.
    \]
    Then we have   
    \begin{align}\label{eq pr conv L2}
        \big(\E {\cal W}^2_2(\hn_t, \hm)\big)^\frac12 = O\left( t^{-\zeta_a} \right).
    \end{align}
    with $\zeta_a = \frac{a}{2(2(d+2)+a(d+3))}$. 
    \color{jfc}
    If $\Sigma$ is uniformly elliptic then one can set $\zeta_a = \frac{a}{2(d+2)(a+2)}$.
    \\
    \color{black}
 \noindent $(b)$    If moreover, we assume that 
    \begin{align} \label{eq ass exp moment L2}
        \exists\, \lambda>0\; \mbox{ such that }\; \int e^{\lambda |y|^2}  \hm(dy)<+\infty ,  
    \end{align}
    (e.g. because $\Sigma$ is bounded, see Lemma~\ref{le sigma bounded exp moments}) then 
    \begin{align}\label{eq pr conv L2 exp moment}
        \big(\E\, \cW_{2}^2(\hn_t,\hm)\big)^\frac12 = O \left(t^{- \zeta_a} 
        (\log t)^{\kappa_a}\right).
    \end{align}
    \color{jfc} with 
    $(\zeta_a,\kappa_a) = \left(\frac1{2(d+3)},\frac{d+2}{2(d+3)}\right)$. If $\Sigma$ is uniformly elliptic then one can set $(\zeta_a,\kappa_a) = \left(\frac1{2(d+2)},\frac12\right)$.
    \color{black}
\end{Proposition}

\smallskip
\color{black}
\proof
\noindent {\sc Step~1} ({\em Decomposition}).
% \color{yellow}
% Let $\rho $ be a positive mollifier with support   $[-1,1]^d$ and let $\zeta$ be a   $\rho(y)\lambda_d(dy)$-distributed random vector. Let $\ve>0$ and let  us define   
%  $$
%  \rho_{\ve}(u) = \ve^{-d}\rho\Big(\frac{u}{\ve}\Big).
%  $$
% Note that $\Lip{\rho_{\ve}} = \ve^{-(d+1)}\Lip{\rho}=  \ve^{-(d+1)}\|\nabla \rho\|_{\sup}$. 
%
% \noindent To control ${\cal W}_2(\hn_t, \hm)$ we first classically smoothen the problem through the triangle inequality as follows
% \begin{align*}
% {\cal W}_2(\hn_t, \hm) &\le {\cal W}_2(\hn_t, \hn_t*[\ve \zeta]) +  {\cal W}_2(\hm_t, \hm *[\ve\zeta])+{\cal W}_2(\hn_t *[\ve\zeta], \hm *[\ve\zeta]) .
% \end{align*}
% One checks that  the distributions $\hn_t*[\ve \zeta]$ and $ \hm *[\ve\zeta]$ have densities respectively  given for every $t\ge 0$ by 
% \[
% \frac{d\hn_t *[\ve\zeta]}{d\lambda_d} =   \frac{1}{t}\int_0^t\rho_{\ve}(\xi-\whX_s) \ud s \quad\mbox{ and }\quad  
% \frac{d\hm_t *[\ve\zeta]}{d\lambda_d} = \int_{\R^d}\rho_{\ve}(\xi-y)\hm(dy).  
% \]
% If we set 
% \begin{equation}\label{eq:dfnphieps}
% \phi_{\ve}(\xi, x) =\rho_{\ve}(\xi-x)-\int_{\R^d}\rho_{\ve}(\xi-y)\hm(dy) 
% \end{equation}
% then, it follows from~\eqref{eq:ident2 new} that  
% \color{black}
From Lemma \ref{le smoothing generic}, taking expectation in \eqref{eq:decomp0 le}, we get
\begin{align}\label{eq:decomp0}
%\nonumber 
\E\, {\cal W}_2(\hn_t, \hm)^2
%&\le 6\ve^2 \|\zeta\|^2_2 + 9\esp{\int_{\R^d} |\xi|^2  \Big|\frac1t\int_0^t \phi_{\ve} (\xi, \whX_s) ds \Big| \ud \xi}  \\
& \le 6\ve^2 \|\zeta\|^2_2 +9\int_{\R^d} |\xi|^2  \esp{\Big|\frac1t\int_0^t \phi_{\ve} (\xi, \whX_s) ds \Big|} \ud \xi .
\end{align}
We split the  above integral into two terms as follows: for every $R \ge 2\,\ve \, \sqrt{d}$,
\begin{equation}\label{eq:decomp1}
 \int_{\R^d} |\xi|^2\esp{ \Big|\frac1t\int_0^t \phi_{\ve} (\xi, \whX_s) ds \Big| }\ud \xi  = I_1(R)  +I_2(R) 
\end{equation} 
 with
 \[
I_1(R)=  \int_{\{|\xi| \le R\}} |\xi|^2\esp{\Big|\frac1t\int_0^t \phi_{\ve} (\xi, \whX_s) ds \Big|} \ud \xi\; \mbox{ and }\; I_2(R)= \int_{\{|\xi|> R\}} |\xi|^2 \esp{\Big|\frac1t\int_0^t \phi_{\ve} (\xi, \whX_s) ds \Big|} \ud \xi .
\]

\noindent {\sc Step 2.1} ({\em Term $I_1(R)$}). It is clear that
\[
I_1(R)\le R^{d+2} V_d \sup_{\xi\in \R^d} \esp{\Big|\frac1t\int_0^t \phi_{\ve} (\xi, \whX_s) ds \Big|}.
\]
where $V_d = \lambda_d\big( B_d(0,1)\big)$ (hyper-volume of the Euclidean unit ball). 

\noindent Applying then Lemma~\ref{le only lip but no ps} {\color{jfc}(we do not assume that $\Sigma$ is uniformly elliptic for now)}, 
there exists a real constant $C= C_{\rho, B,\Sigma,\a}>0$ such that
\begin{align}\label{eq I1 esp}
    I_1(R)\le C\, [\rho_{\ve}]_{\rm Lip}\frac{ R^{d+2}}{\sqrt{t}} \le C\, \frac{ R^{d+2}}{\sqrt{t}\,\ve^{d+1}}.
\end{align}

\noindent {\sc Step~2.2} ({\em Term $I_2(R)$}). We first note, using Fubini's Theorem at each line  and the elementary inequality $(x+y)^2 \le 2(x^2+y^2)$, $x,y\ge 0$, in the third line that 
\begin{align*}
I_2(R) &\le \int_{\{|\xi|>R\}}|\xi|^2\bigg( \esp{\frac 1t \int_0^t\rho_{\ve}(\xi-\whX_s)ds} + \int \rho_{\ve}(\xi-y)\hm(dy)\bigg) \ud \xi\\
&\le 2\int_{\{|\xi|>R\}}\esp{|\xi|^2\frac1t \int_0^t\rho_{\ve}(\xi-\whX_s)ds}\ud \xi
\\
&\le 
 4 \Big(\,I_{21}(R) + I_{22}(R) \,\Big)\,,
\end{align*}
with 
\[ 
I_{21}(R)  = \frac1t \int_0^t \int_{\{|\xi|>R\}} \esp{|\xi-\whX_s|^2  \rho_{\ve}(\xi-\whX_s)}\ud \xi ds
\]
 and  
 \begin{align}
    I_{22}(R)  = \frac1t \int_0^t\int_{\{|\xi|>R\}} \esp{|\whX_s|^2\rho_{\ve}(\xi-\whX_s)}\ud \xi ds. \label{eq I22 esp} 
 \end{align}    

Now
\begin{align*}
%\frac1t\int_0^t \int_{\{|\xi|>R\}} \esp{|\xi-\whX_s|^2  \rho_{\ve}(\xi-\whX_s)}\ud \xi ds
I_{21}(R) &\le \frac1t\int_0^t \int_{\R^d} \esp{|\xi-\whX_s|^2  \rho_{\ve}(\xi-\whX_s)}\ud \xi ds
%\\ & 
= \int_{\R^d} |z|^2\rho_{\ve}(z)dz = \ve^{2}\|\zeta\|_2^2,
\end{align*}
where we used the change of variable $\xi = z+\whX_s$ in the inner integral. On the other hand, as $\rho_{\ve}$ is $[-\ve,\ve]^d$-supported and $|\cdot|\le \sqrt{d}|\cdot|_{{\infty}}$,
\begin{align*}
    I_{22}(R) &\le  \int_{\{|\xi|>R\}}\esp{\mbox{\bf1}_{\{|\whX_s-\xi|\le\sqrt{d}\, \ve\}} |\whX_s|^2\rho_{\ve}(\xi-\whX_s)}\ud \xi.
\end{align*}
Note that $\{|\whX_s-\xi|\le\sqrt{d}\, \ve\}\subset \{|\whX_s|\ge R-\ve\, \sqrt{d}\}$ when $|\xi| > R$. Hence, it follows from Markov inequality that 
\begin{align*}
 \int_{\{|\xi|>R\}}\!\!\esp{\mbox{\bf1}_{\{|\whX_s-\xi|\le\sqrt{d}\, \ve\}} |\whX_s|^2\rho_{\ve}(\xi-\whX_s)}\ud \xi&\le \frac{1}{(R-\ve\,\sqrt{d})^a} \int_{\{|\xi|>R\}}\!\!\esp{|\whX_s|^{2+a}\rho_{\ve}(\xi-\whX_s)}\ud \xi\\
 & \le  \frac{1}{(R-\ve\,\sqrt{d})^a}\esp{|\whX_s|^{2+a}\int_{\{|\xi|>R\}}\!\!\rho_{\ve}(\xi-\whX_s)\ud \xi}\\
 &\le   \frac{1}{(R-\ve\,\sqrt{d})^a}\esp{|\whX_s|^{2+a}}\int_{\R^d}\rho_{\ve}(\xi)\ud \xi\\
 &= \frac{1}{(R-\ve\,\sqrt{d})^a} \esp{|\whX_s|^{2+a}}\\
 &\le \frac{2^a}{R^a}\int |y|^{2+a}\hm(dy) \quad \mbox{since $\,\ve \sqrt{d} \le R/2$.}
\end{align*}
This shows that, there exists a real constant $C= C_{B,\Sigma,\rho,a}\!\in (0,+\infty)$ such that 
\[
I_2(R)\le C(\ve^2 + R^{-a}).
\]

\noindent {\sc Step~2.3} ({\em Conclusion}). Plugging the above bounds for $I_i(R)$, $i=1,2$, into~\eqref{eq:decomp1} and then in ~\eqref{eq:decomp0} yields, for $t \ge 1$, 
\[
 \esp{ {\cal W}_2(\hn_t, \hm)^2}\le C\big(\ve^2 + R^{-a} + R^{d+2} t^{-1/2}\ve^{-(d+1)}\big)
\]
for some real constant $C=  C_{B,\Sigma,\rho,a}\!\in (0,+\infty)$. For $t \ge 1$, assume $R$ and $\ve$ are of the form 
\[
R= 2\sqrt{d}\,t^{\frac{1}{2(d+2)+a(d+3)}}\quad \mbox{ and }\quad \ve = t^{-\frac{a}{2(2(d+2)+a(d+3))}}
\]
so that $R \ge 2 \ve \sqrt{d}$, then
\[
 \esp{ {\cal W}^2_2(\hn_t, \hm)}^\frac12 \le C_{a,d}\, t^{-\frac{a}{2(2(d+2)+a(d+3))}}.
\]
 
%Note that, as $a\to +\infty$, $\frac{a}{2(2(d+2)+a(d+3))}\to \frac{1}{2(d+3)}$.

\color{jfc}
\noindent {\sc Step~2.4} ({\em Case of uniformly elliptic $\Sigma$}). The main change is in equation \eqref{eq I1 esp}. Indeed,
applying then Lemma~\ref{le lip wf unif ellip bounded} when $\Sigma$ is uniformly elliptic, we have that 
there exists a real constant $C= C_{\rho, B,\Sigma,\a}>0$ such that
\begin{align}\label{eq I1 esp new}
    I_1(R)\le C\, \|\rho_{\ve}\|_{\sup}\frac{ R^{d+2}}{\sqrt{t}} \le C\, \frac{ R^{d+2}}{\sqrt{t}\,\ve^{d}}.
\end{align}
Following the same stream of computations as above, we then obtain as conclusion:
for $t \ge 1$, 
\[
 \esp{ {\cal W}_2(\hn_t, \hm)^2}\le C\big(\ve^2 + R^{-a} + R^{d+2} t^{-1/2}\ve^{-d}\big)
\]
for some real constant $C=  C_{B,\Sigma,\rho,a}\!\in (0,+\infty)$. For $t \ge 1$, assume now that $R$ and $\ve$ are of the form 
\[
R= 2\sqrt{d}\,t^{\frac{1}{(d+2)(a+2)}}\quad \mbox{ and }\quad \ve = t^{-\frac{a}{2(d+2)(a+2)}}
\]
so that $R \ge 2 \ve \sqrt{d}$, then
\[
 \esp{ {\cal W}^2_2(\hn_t, \hm)}^\frac12 \le C_{a,d}\, t^{-\frac{a}{2(d+2)(a+2)}}.
\]
\color{black}
\medskip
\noindent {\sc Step~3.1} We now turn to the case where the stationary measure has exponential moments, namely~\eqref{eq ass exp moment L2} holds. In this case, we follow the proof above but we handle the term $I_{22}(R)$ slightly differently to take advantage  of the exponential moment. Namely, according to~\eqref{eq I22 esp}, we compute 
\begin{align*}
    I_{22}(R) &\le  \int_{\{|\xi|>R\}}\esp{\mbox{\bf1}_{\{|\whX_s|>R - \ve \sqrt{d}\}} |\whX_s|^2\rho_{\ve}(\xi-\whX_s)}\ud \xi
\end{align*}
since $\rho_{\ve}$ is $[-\ve,\ve]^d$-supported. By assumption, we have that $R - \ve \sqrt{d} \ge \frac{R}2$. Hence, it follows  successively from Markov inequality and~\eqref{eq ass exp moment L2} that
\begin{align}
  \nonumber   I_{22}(R)  \le 
    \esp{e^{\frac{\lambda}2|\whX_s|^2}|\whX_s|^2e^{-\frac{\lambda R^2}4}\int_{\{|\xi|>R\}}\rho_{\ve}(\xi-\whX_s) \ud \xi }&\le C_\lambda e^{-\frac{\lambda R^2}4} \esp{e^{\lambda|\whX_s|^2}}
    \\
    &\le C_\lambda e^{-\frac{\lambda R^2}4}.
\end{align}
 In this case, we thus obtain that 
\begin{align*}
    I_2(R) \le C(\ve^2 + e^{-\frac{\lambda R^2}4}),
\end{align*}
which combined with the bound~\eqref{eq I1 esp} leads to, for $t \ge 1$, 
\begin{align}\label{eq l2 conv expo moment first}
    \esp{ {\cal W}_2(\hn_t, \hm)^2}\le C\big(\ve^2 + e^{-\frac{\lambda}4R^2} + R^{d+2} t^{-1/2}\ve^{-(d+1)}\big).
\end{align}
We then first set $\ve = t^{-\frac{1}{2(d+3)}}R^{\frac{d+2}{d+3}}$ to get
\[
    \esp{ {\cal W}_2(\hn_t, \hm)^2}\le C\big( 
        t^{-\frac{1}{d+3}}R^{2\frac{d+2}{d+3}} + e^{-\frac{\lambda}4R^2}  
    \big).
    \]
    Then, choosing $R=\frac4{\lambda\wedge 1}\log\left(t^{1/(d+3)}\right)^\frac12$, we observe that for $t$ large enough $R \ge 2 \ve \sqrt{d}$ and then we obtain 
    \[
        \esp{ {\cal W}_2(\hn_t, \hm)^2}=O\left(\frac{(\log t)^{\frac{d+2}{d+3}}}{t^{1/(d+3)}} \right), 
        \]
      \color{jfc}  which concludes the proof for the case when $\Sigma$ is not uniformly elliptic.
\\
\noindent {\sc Step~3.2} When $\Sigma$ is uniformly elliptic, relying on \eqref{eq I1 esp new} instead of \eqref{eq I1 esp}, we obtain  
\begin{align}\label{eq l2 conv expo moment second}
    \esp{ {\cal W}_2(\hn_t, \hm)^2}\le C\big(\ve^2 + e^{-\frac{\lambda}4R^2} + R^{d+2} t^{-1/2}\ve^{-d}\big).
\end{align}
We then set $\ve = t^{-\frac{1}{2(d+2)}}R$ and $R=\frac{4}{\lambda \wedge 1}\log(t^{\frac1{d+2}})^\frac12$ to get
\[
        \esp{ {\cal W}_2(\hn_t, \hm)^2}=O\left(\frac{\log t}{t^{1/(d+2)}} \right),
        \]
which concludes the proof.
\eproof
\color{black}

\section{An $a.s.$  convergence rate}

\label{se as conv}

\color{black}
In this section, we prove almost sure rates for the convergence of the empirical measure towards the stationary distribution in Wasserstein distance. As in the previous section, we rely on the novel use of coboundary functions for a specific class of functions. However here, in this almost sure convergence setting, we have to pay special attention to martingale terms. Among other things, we call upon some concentration inequalities and stochastic Kronecker Lemmas presented in the appendix. 

\color{black}
\begin{Proposition}\label{pr W2 as conv} Let \HYP{L} and \HYP{C} hold and
    assume that $\Sigma$ is uniformly elliptic. Let $a>0$ from Lemma \ref{le sigma bounded exp moments} such that
    \begin{align}\label{eq moment ass as}
         \int|y|^{2+a}  \hm(dy)<+\infty.    
    \end{align}
    Then 
    \begin{align}\label{eq pr conv as}
        \cW_2(\hn_t,\hm) = o_\eta\big(t^{-\zeta_a}\log^{\frac12+\eta}(t)\big)\,, \quad \mbox{ for every $\eta > 0$},   
    \end{align}
    where %$\zeta_a = \frac{a^2}{2(2+a) \{2(d+2) + (d+3)(a+d) \}}$.
    \textcolor{myred}{$ \zeta_a =\frac{a^2}
{2(a+2)\left(d(d+3)+(a+2)(d+2)\right)}$}.
\end{Proposition}

\color{black}
\proof 
\noindent {\sc Step~1} ({\em Decomposition}). 
% Let $\rho $ be a positive mollifier with support   $[-1,1]^d$ and let $\zeta$ be a   $\rho(y)\lambda_d(dy)$-distributed random vector. Let $\ve>0$ and let  us define   
%  $$
%  \rho_{\ve}(u) = \ve^{-d}\rho\Big(\frac{u}{\ve}\Big).
%  $$
% Note that $[\rho_{\ve}]_{\rm Lip} = \ve^{-(d+1)}[\rho]_{\rm Lip}=  \ve^{-(d+1)}\|\nabla \rho\|_{\sup}$ and $\|\rho_{\ve} \|_{\sup} \le \ve^{-d}\|\rho \|_{\sup}$. 
%
%\noindent 
We specify a sequence of smoothing parameters $(\varepsilon_n)_{n \ge 1}$, namely 
\begin{align}\label{eq de epsilon unbounded sigma}
\varepsilon_n = n^{-r}
\end{align}
for some $r>0$ to be fixed later on.\\
From now on, we work assuming that $t \ge 1$. Invoking Lemma \ref{le smoothing generic} (with the above smoothing parameter sequence), we have from \eqref{eq:decomp0 le} that
% To control ${\cal W}_2(\hn_t, \hm)$, we first classically smoothen the problem through the triangle inequality as follows
% \begin{align}\label{eq decomp Wasserstein unbounded}
% {\cal W}_2(\hn_t, \hm) &\le {\cal W}_2(\hn_t, \hn_t*[\ve_{\ent{t}} \zeta]) +  {\cal W}_2(\hm, \hm *[\ve_{\ent{t}}\zeta])+{\cal W}_2(\hn_t *[\ve_{\ent{t}}\zeta], \hm *[\ve_{\ent{t}}\zeta]) .
% \end{align}
% \noindent One checks that  the distributions $\hn_t*[\ve_{\ent{t}} \zeta]$ and $ \hm *[\ve_{\ent{t}}\zeta]$ have densities respectively  given  by 
% \[
% \frac{d\hn_t *[\ve_{\ent{t}}\zeta]}{d\lambda_d}(\xi) =   \frac{1}{t}\int_0^t\rho_{\ve_{\ent{t}}}(\xi-\whX_s) ds \quad \mbox{
% and }\quad 
% \frac{d\hm_t *[\ve_{\ent{t}}\zeta]}{d\lambda_d}(\xi) = \int_{\R^d}\rho_{\ve_{\ent{t}}}(\xi-y)\hm(dy).  
% \]
% If we set 
% \begin{equation}\label{eq:dfnphieps as}
% \phi_{\ve_{\ent{t}}}(\xi, x) =\rho_{\ve_{\ent{t}}}(\xi-x)-\int_{\R^d}\rho_{\ve_{\ent{t}}}(\xi-y)\hm(dy) 
% \end{equation}
% then, it follows from~\eqref{eq decomp Wasserstein unbounded} and \eqref{eq:ident2 new} 
\begin{align}\label{eq:decomp0 as unbounded}
 {\cal W}_2(\hn_t, \hm)^2
&\le 6\ve_{\ent{t}}^2 \|\zeta\|^2_2 + 9 \int_{\R^d} |\xi|^2  \Big|\frac1t\int_0^t \phi_{\ve_{\ent{t}}} (\xi, \whX_s) ds \Big| \ud \xi.   
%& \le 6\ve^2 \|\zeta\|^2_2 +3\bigg(\int_{\R^d} |\xi|^2  \E\, \Big|\frac1t\int_0^t \phi_{\ve} (\xi, X_s) ds \Big| \ud \xi \bigg),
\end{align}
%where we used $\|\cdot\|_1\le \|\cdot\|_2$ and Fubini's theorem in the second line.
%
%\smallskip 
We now consider an increasing sequence $(R_n)_{n \ge 0}$ of truncation parameters with 
\begin{align}\label{eq de R unbounded}
    R_n = 2\sqrt{d}\,n^{\kappa}, \text{ for } n \ge 1 \text{ and } R_0 := R_1, 
\end{align}
for some $\kappa>0$ to be fixed later on.\\
We observe
\begin{align}\label{eq small truc pour plus tard bis unbounded}
    R_n \ge 2 \sqrt{d}\,\varepsilon_{n'} \quad \text{ for } \quad n \ge 0,n'\ge 1.
\end{align}
This will prove useful at some point in the proof below. 

%  \textcolor{magenta}{\bf verifier pourquoi on part de $2\epsilon \sqrt{d}$ pour la dernière partie de la preuve...}:
% \begin{align*}
%     \R_+ \ni t \mapsto R(t) \in (2\epsilon \sqrt{d},+\infty),
% \end{align*}
% such that $R(t)\sim t^\kappa$, for some $\kappa > 0$ to be chosen later on.
% \textcolor{red}{bien reprendre parametrisation...}
\noindent Now, we split the  above integral into two terms using the sequence $(R_n)_n$ as follows: 
\begin{equation}\label{eq:decomp1 as unbounded}
 \int_{\R^d} |\xi|^2\E\, \Big|\frac1t\int_0^t \phi_{\ve_{\ent{t}}} (\xi, \whX_s) ds \Big| \ud \xi  = I_1(t)  +I_2(t) 
\end{equation} 
 with
 \[
I_1(t)=  \int_{\{|\xi| \le R_{\ent{t}}\}} |\xi|^2\Big|\frac1t\int_0^t \phi_{\ve_{\ent{t}}} (\xi, \whX_s) ds \Big| \ud \xi\; \mbox{ and }\; I_2(t)= \int_{\{|\xi|> R_{\ent{t}}\}} |\xi|^2 \Big|\frac1t\int_0^t \phi_{\ve_{\ent{t}}} (\xi, \whX_s) ds \Big| \ud \xi .
\]

\noindent {\sc Step 2} ({\em Term $I_1(t)$}). 
%\textcolor{yellow}{I do the computation with $d=1$ to start with.} 
% Let $(\gamma_n)$ a decreasing sequence of steps and $\Gamma_n = \sum_{k=1}%^n\gamma_k$, it is assumed that $\Gamma_n \rightarrow +\infty$ when $n \rightarrow %+\infty$. Recall notation: $\t$, $N(t)$.
% \textcolor{red}{on peut prendre un pas de temps constant !! $\Gamma_n \sim n$ %simplement donc...}
%
We first observe that 
\begin{align}
   \nonumber  I_1(t)&\le \int_{|\xi|\le R_{\ent{t}}}|\xi|^2\Big |\frac1{t}\int_{\ent{t}}^t\phi_{\ve_{\ent{t}}}(\xi,\whX_s)\ud s\Big|\ud \xi + \int_{|\xi|\le R_{\ent{t}}}|\xi|^2 \Big|\frac1{\ent{t}}\int_{0}^{\ent{t}}\phi_{\ve_{\ent{t}}}(\xi,\whX_s)\ud s\Big|\ud \xi
    \\
    &\le C\frac{R_{\ent{t}}^{d+2}}{t\ve_{\ent{t}}^{d}}+ \tilde{I}_1(\ent{t})
    % \int_{|\xi|\le R(\Gamma_{N(t)+1})}|\xi|^2 |\frac1{\Gamma_{N(t)}}\int_{0}^{\Gamma_{N(t)}}\phi_{\ve(t)}(\xi,X_s)\ud s|\ud \xi
    % =:I_{12}(t)
\end{align} 
where, for $n \ge 1$, 
\[
    \tilde{I}_1(n) := \int_{|\xi|\le R_n}|\xi|^2 \Big|\frac1{n}\int_{0}^{n}\phi_{\ve_n}(\xi,\whX_s)\ud s\Big|\ud \xi\,.
    \]
We divide $[-R_n,R_n]$ into $K_n$ interval of size $\delta_n := 2R_n/K_n$ and subsequently $[-R_n,R_n]^d$ is covered with $K_n^d$ hypercube of length $\delta_n$. For a given indexation $k=0,\dots,K_n^d-1$, we denote these hypercubes $A^n_k$.
Then,
\begin{align}
    \tilde{I}_{1}(n) \le R_n^2 \sum_{k=0}^{K_{n}^d-1}\int_{\xi \in A^n_k} \Big|\frac1{n}\int_{0}^{n}\phi_{\ve_n}(\xi,\whX_s)\ud s\Big|\ud \xi.
\end{align}
\color{myred}As $\phi_{\ve_n}$ is $\ve_n^{-(d+1)}[\rho]_{\rm Lip}$-Lipschitz, we observe that  for every $\xi \in A^n_k$, 
\begin{align}
    |\phi_{\ve_n}(\xi,\whX_s)-\phi_{\ve_n}(\xi^n_k,\whX_s)| \le C\frac{\delta_n}{\ve_n^{d+1}},
\end{align}
\color{black}
where $\xi^n_k$ is the center point of $A^n_k$. Since $\int_{\xi \in A^n_k}\ud \xi \le C\delta_n^d$,
%(the lower corner of $A^n_k$ say).
%This leads to, observing 
it follows that 
\begin{align}
    \tilde{I}_{1}(n) \le C\bigg(\frac{R_n^{d+2}}{\textcolor{myred}{\ve_n^{d+1}}} \delta_n + R_n^2 \delta_n^d \sum_{k=0}^{K_{n}^d-1} \Big|\frac1{n}\int_{0}^{n}\phi_{\ve_n}(\xi^n_k,\whX_s)\ud s\Big|\bigg)
\end{align}
and 
 \begin{align}
    \tilde{I}_{1}(n) \le C \bigg(\frac{R_n^{d+2}}{\textcolor{myred}{\ve_n^{d+1}}} \delta_n + R_n^{d+2} \sup_{k=0,\ldots,K_n^d-1}\Big|\frac1{n}\int_{0}^{n}\phi_{\ve_n}(\xi^n_k,\whX_s)\ud s\Big|\bigg).
 \end{align}
Now, invoking Proposition~\ref{le unif elliptic cobord} and Lemma~\ref{le lip wf unif ellip bounded}, we have (writing $\omega_{\ve_n}$ for $\omega_{\phi_{\ve_n}}$)
\begin{align}
   \nonumber \Big |\frac1{n}\int_{0}^{n}\phi_{\ve_n}(\xi^n_k,\whX_s)\ud s\Big|
    &\le \frac1{n}|\omega_{\ve_n}(\xi^n_k,\whX_{n})-\omega_{\ve_n}(\xi^n_k,\whX_0)\Big|
    + \frac1{n}\Big|\int_{0}^{n}\partial_x \omega_{\ve_n}(\xi^n_k,\whX_s)\Sigma(\whX_s)\ud W_s\Big|
    \\
    &\le  C\frac{|\whX_{n}-\whX_{0}|}{\textcolor{jfc}{\ve_n^{d}}n}
    + \frac{\psi(n)}{n}\Big|\frac1{\psi(n)}\int_{0}^{n}\partial_x \omega_{\ve_n}(\xi^n_k,\whX_s)\Sigma(\whX_s)\ud W_s\Big|
\end{align}
where $\psi(n)=n^{\theta}(\log n)^{\frac{1+\eta}{a+2}}$, $n\ge 2$,  for some $ \theta \!\in [\frac12 , 1)$ to be chosen later on and $\eta>0$.
Hence
\begin{align*}
%\label{eq temp Itilde un}
    \tilde{I}_{1}(n) \le C \frac{R_n^{d+2}}{\textcolor{jfc}{\ve_n^{d}}} \left(\frac{\delta_n}{\textcolor{myred}{\ve_n}}+\frac{|\whX_{n}-\whX_{0}|}{n} +  \frac{\psi(n)}{{n}}\sup_{k=0,\ldots,K_n^d-1}\Big|\frac1{\psi(n)}\int_{0}^{n}\textcolor{jfc}{\ve_n^{d}}\partial_x \omega_{\ve_n}(\xi^n_k,\whX_s)\Sigma(\whX_s)\ud W_s\Big|\right).
\end{align*}
% \textcolor{red}{Je me suis arrêté là pour cette étape. En effet, quand on passe de~\eqref{eq temp Itilde un} à~\eqref{eq majo violente} il y a $\frac1{\ve_n^{d+1}}$ qui est mis en factor et donc pour le terme martingale ci-dessus on devrait avoir $\ve_n^{d+1}$ qui compense $|\partial_x \omega_{\ve_n}|_\infty$. Du coup j'ai eu l'impression qu'on avbait un peu plus de place pour obtenir une vitesse de convergence (toute pourrie) dans le cas sigma non-bornée... je me suis endormi en essayant d'optimiser les paramètres... je n'ai donc pas fini de vérifier l'étape 2 et donc je n'ai pas repris l'étape 4 de conclusion.}
Now we want to prove for an appropriate choice of $\theta$ that 
\begin{align}\label{eq as conv gyorfi}
    \sup_{k=0,\ldots,K_n^d-1}\Big|\frac1{\psi(n)}\int_{0}^{{n}}\textcolor{jfc}{\ve_n^{d}}_n\partial_x \omega_{\ve_n}(\xi^n_k,\whX_s)\Sigma(\whX_s)\ud W_s\Big|\rightarrow 0 \text{ as } n \rightarrow +\infty.
\end{align}
Let $\ell > 0$. Denote for every $k\!\in \{0,\ldots,K_n^d-1\}$, 
\[
Z^{n,k}_s =     \textcolor{jfc}{\ve_n^{d}}\partial_x \omega_{\ve_n}(\xi^n_k,\whX_s)\Sigma(\whX_s).
\]
Since $|Z^{n,k}_s|\le C\textcolor{jfc}{\|\rho\|_{\sup}}(1+|\whX_s|)$, we deduce from~\eqref{eq moment ass as} that 
%\begin{align}\label{eq control integral}
    $\esp{|Z^{n,k}_s|^{2+a}} \le C$. 
%\end{align}
Applying~\eqref{eq encore merci bdg} in Lemma~\ref{le hoeffding like}$(b)$, we obtain
%denoting $q = \frac{a}2+1$,
\begin{align}
    \P\left(\Big|\frac1{\psi(n)}\int_{0}^{n}Z^{n,k}_s\ud W_s\Big| > \ell\right)
    \le  \frac{C}{\ell^{a+2}n^{(a+2)(\theta-\frac12)}( \log n)^{1+\eta}}.
\end{align}
Then,
\begin{align}
    \P\left(\frac1{\psi(n)}\sup_{k=0,\ldots,K_n^d-1}\Big|\int_{0}^{{n}}Z^{n,k}_s\ud W_s\Big| > \ell\right)
    \le \frac{C K_n^d}{\ell^{a+2}n^{(a+2)(\theta-\frac12)}( \log n)^{1+\eta}}.
\end{align}
We now set $K_n = \lfloor n^{p}\rfloor\sim n^p$ for some $p>0$ to be determined later on.
%
%     Set $\ve_n = \lfloor n^{-r}\rfloor\sim n^{-r}$ and,  from the previous discussion in the $L^2$-convergence case, set 
%   \begin{align}
%   \theta = \frac12 +r(d+1).
%   \end{align}
%   Then $\ve_n^{2(d+1)} n^{2\theta-1} \sim 1$. 
  Hence, 
  %as soon as $K_n \sim n^p$ for some $p$,
% and $\ve_n^{2(d+1)} n^{2\theta-1} \sim n^{\eta}$ for some $\eta>0$, 
we obtain
\begin{align}\label{eq majo sans scrupule unbounded}
 \hskip-0,3cm    \sum_{n=1}^{+\infty}\P\left(\frac1{\psi(n)}\sup_{k=0,\ldots,K_n^d-1}\Big|\int_{0}^{n}Z^{n,k}_s\ud W_s\Big| > \ell\right)
    \le \frac{C}{\ell^{a+2}}\sum_{n=1}^{+\infty}   \frac{1}{n^{(a+2)(\theta-\frac12)-pd}( \log n)^{1+\eta}}.
    % < +\infty
\end{align}
We thus impose $0<p< \frac{a}{2d}$ and set  
\begin{align} \label{eq de theta p}
    \theta=\frac12  + \frac{1+pd}{a+2}\!\in (\tfrac 12, 1).
\end{align}
%    imposing that $\theta < 1$ i.e. 
%    \begin{align} \label{constraint through theta}
%       0 \le p < \frac{a}{2d}\;.
%    \end{align}
    Then, we deduce from~\eqref{eq majo sans scrupule unbounded},
    \begin{align}
        \sum_{n=1}^{+\infty}\P\left(\frac1{\psi(n)}\sup_{k=0,\ldots,K_n^d-1}\Big|\int_{0}^{n}Z^{n,k}_s\ud W_s\Big| > \ell\right)
        \le \frac{C}{\ell^{a+2}}\sum_{n=1}^{+\infty}   \frac{1}{n( \log n)^{1+\eta}} < +\infty\,.
    \end{align}
     It follows that, for every (rational) $\ell> 0$, the series in the r.h.s. of the above equation is converging so that, by the Borel-Cantelli Lemma, 
\begin{align}\label{eq you deserve it0}
    \limsup_n \frac1{\psi(n)}\sup_{k=0,\ldots,K_n^d-1}\Big|\int_{0}^{{n}}Z^{n,k}_s\ud W_s\Big| \le \ell \quad \P\mbox{-} a.s.
\end{align}
or, equivalently,
\begin{align}\label{eq you deserve it}
    \limsup_n\frac1{\psi(n)} \sup_{k=0,\ldots,K_n^d-1} \Big|\int_{0}^{{n}}Z^{n,k}_s \ud W_s\Big| = 0 \quad \P\mbox{-} a.s.
\end{align}

Furthermore, we know from Lemma~\ref{le control SDE}$(b)$ setting $a'=a\wedge 2$ that
\begin{align}
    \frac{|X_{n}-X_{0}|}{n} = o_{\eta}\big(n^{-\frac{1+a'}{2+a'}}(\log n)^{\frac{1}{1+a'/2}+\eta}\big) =o_{\eta}\big(n^{-\frac12}(\log n)^{{1+\eta}}\big) =o\Big(\frac{\psi(n)}{n}\Big)
    \end{align}
%and this term can be neglected in regards to $\frac{\psi(n)}{n}$, according to the definition of $\theta$ in~\eqref{eq de theta p}.
since $\theta -1>-\frac 12$. This leads to
\begin{align}\label{eq majo violente}
    \tilde{I}_{1}(n) \le C \frac{R_n^{d+2}}{\textcolor{jfc}{\ve_n^{d}}} \left(\frac{\delta_n}{\textcolor{myred}{\ve_n}}
    +  o\Big(\frac{\psi(n)}{n}\big) 
    \right) .
\end{align}

Since $R_n = n^\kappa$ then $\delta_n=\frac{2R_n}{K_n} \sim 2 n^{\kappa-p}$ and
we obtain 
\begin{align}\label{eq majo almost done}
    \tilde{I}_{1}(n) \le C 
    \left(n^{(d+3)\kappa + \textcolor{myred}{r(d+1)}-p} + 
    n^{(d+2)\kappa + \textcolor{jfc}{rd}+\theta-1}o(\log n^{\frac1{a+2} + \eta}) 
    \right) .
\end{align}
Finally, one  checks  that the term  $\frac{R_{\ent{t}}^{d+2}}{t\ve_{\ent{t}}^{d}}= O(\ent{t}^{(d+2)\kappa+ \textcolor{jfc}{rd}+\theta-1 -(r+\theta)})= o \big( \tilde{I}_{1}(\ent{t})\big) $.  
%and we still have~\eqref{eq you deserve it}.
So the conclusion for this step is  that, for every $\eta > 0$,
\begin{align}\label{eq conclu term 1}
    I_1(t) \le C \left(t^{(d+3)\kappa + \textcolor{myred}{r(d+1)}-p} + 
    t^{(d+2)\kappa + \textcolor{jfc}{rd}+\theta-1}o(\log(t)^{\frac1{a+2} + \eta}) 
    \right).
\end{align}

\medskip
\noindent {\sc Step~3} ({\em Term $I_2$}). 
We remark that
%, using Fubini's Theorem at each line  and the elementary inequality $(a+b)^2 \le 2(a^2+b^2)$, $a,b\ge 0$, in the third line that 
\begin{align}\label{eq:decompI_2}
% I_2(R(t)) &\le \int_{\{|\xi|>R(t)\}}|\xi|^2\bigg( \frac 1t \int_0^t\rho_{\ve}(\xi-X_s)ds + \int \rho_{\ve}(\xi-y)\hm(dy)\bigg) \ud \xi\\
\nonumber I_2(t) &\le 
\int_{\{|\xi|>R_{\ent{t}}\}}\frac{|\xi|^2}t \int_0^t\rho_{\ve_{\ent{t}}}(\xi-\whX_s)\ud s \ud \xi 
+ \int_{\{|\xi|>R_{\ent{t}}\}}|\xi|^2\int \rho_{\ve_{\ent{t}}}(\xi-y)\hm(\d y)\ud \xi\\
&=:I_{21}(t)+I_{22}(t).
% &\le 2\int_{\{|\xi|>R(t)\}} |\xi|^2 \bigg(\frac1t \int_0^t\rho_{\ve}(\xi-X_s)ds\bigg)\ud \xi\\
% & \le  4\bigg[\frac1t\int_0^t \int_{\{|\xi|>R\}} \E\, |\xi-X_s|^2  \rho_{\ve}(\xi-X_s)\ud \xi ds + \frac1t \int_0^t\int_{\{|\xi|>R\}} \E\, |X_s|^2\rho_{\ve}(\xi-X_s)\ud \xi ds  \bigg].
\end{align}
$\blacktriangleright$ We first study $I_{22}$: We compute 
\begin{align*}
    \int_{\{|\xi|>R_{\ent{t}}\}}|\xi|^2\int \rho_{\ve_{\ent{t}}}(\xi-y)\hm(\ud y)\ud \xi
    \le& 2   \int_{\{|\xi|>R_{\ent{t}}\}}\int |\xi-y|^2\rho_{\ve_{\ent{t}}}(\xi-y)\hm(\ud y)\ud \xi 
    \\
    &+ 2 \int_{\{|\xi|>R_{\ent{t}}\}}\int |y|^2\rho_{\ve_{\ent{t}}}(\xi-y)\hm(\ud y)\ud \xi.
\end{align*}
For the first term in the r.h.s. of the previous inequality, we have 
\begin{align*}
    \int_{\{|\xi|>R_{\ent{t}}\}}\int |\xi-y|^2\rho_{\ve_{\ent{t}}}(\xi-y)\hm(\ud y)\ud \xi
    &\le
    \int\int  |\xi-y|^2\rho_{\ve_{\ent{t}}}(\xi-y)\hm(\ud y)\ud \xi
    \\
    &
    = \int|z|^2\rho_{\ve_{\ent{t}}}(z)\ud z = \ve_{\ent{t}}^2 \esp{|\zeta|^2}.
\end{align*}
For the second term,  first note that, as $\rho_{\ve_{\ent{t}}}$ is $[-\ve_{\ent{t}},\ve_{\ent{t}}]^d$-supported and $|\cdot|\le \sqrt{d}|\cdot|_{{\infty}}$,  
\begin{align*}
    \int_{\{|\xi|>R_{\ent{t}}\}}\int |y|^2 \rho_{\ve_{\ent{t}}}(\xi-y)\hm(\ud y)\ud \xi
    \le   
    \int_{\{|\xi|>R_{\ent{t}}\}}\int |y|^2\1_{\set{|y-\xi|<\sqrt{\ve_{\ent{t}}}d}}\rho_{\ve_{\ent{t}}}(\xi-y)\hm(\ud y) \ud \xi.
\end{align*}
We also observe that 
$$
\set{y:|y-\xi|\le \ve_{\ent{t}}\sqrt{d}, |\xi|> R_{\ent{t}}}\subset \set{y:|y|\ge R_{\ent{t}}- \ve_{\ent{t}}\sqrt{d}, |\xi|> R_{\ent{t}}}.$$
It follows then from Markov inequality that 
\begin{align*}
 \int_{\{|\xi|>R_{\ent{t}}\}}\int \mbox{\bf1}_{\{|y-\xi|\le\sqrt{d}\, \ve_{\ent{t}}\}} |y|^2&\rho_{\ve_{\ent{t}}}(\xi-y)\hm(\ud y)\ud \xi 
 \\
 &\le \frac{1}{(R_{\ent{t}}-\ve_{\ent{t}}\,\sqrt{d})^a} \int_{\{|\xi|>R_{\ent{t}}\}} \int_{\R^d}|y|^{2+a}\rho_{\ve}(\xi-y)\hm(\ud y)\,\ud \xi\\
 &\le \frac{1}{(R_{\ent{t}}-\ve_{\ent{t}}\,\sqrt{d})^a}  \int_{\R^d}|y|^{2+a}\int \rho_{\ve_{\ent{t}}}(\xi-y)\ud \xi\,\hm(\ud y)\\
&\le \frac{1}{(R_{\ent{t}}-\ve\,\sqrt{d})^a}  \int_{\R^d}|y|^{2+a}\hm(\ud y)\\
%  & \le  \frac{1}{(R-\ve\,\sqrt{d})^a}\E |X_s|^{2+a}\int_{\{|\xi|>R\}}\rho_{\ve}(\xi-X_s)\ud \xi\\
%  &\le   \frac{1}{(R-\ve\,\sqrt{d})^a}\E |X_s|^{2+a}\int_{\R^d}\rho_{\ve}(\xi)\ud \xi\\
%  &= \frac{1}{(R-\ve\,\sqrt{d})^a} \E |X_s|^{2+a}\\
 &\le \frac{2^a}{R_{\ent{t}}^a}\int |y|^{2+a}\hm(\ud y)
\end{align*}
since  $\ve_{\ent{t}} \sqrt{d} <R_{\ent{t}}/2$, according to~\eqref{eq small truc pour plus tard bis unbounded}. We then conclude for this step that 
\begin{align}\label{eq I22}
    I_{22}(t) \le C\big(\ve^2_{\ent{t}} + R_{\ent{t}}^{-a}\big).
\end{align}
$\blacktriangleright$   We now study $I_{21}$.
We first observe, since $R_n$ is non-decreasing,
\begin{align*}
    \int_{\{|\xi|>R_{\ent{t}}\}}\frac{|\xi|^2}t \int_0^t\rho_{\ve_{\ent{t}}}(\xi-\whX_s)\ud s \ud \xi
    \le &
    \frac{2}t \int_0^t\int_{\{|\xi|>R_{\ent{t}}\}} |\xi - \whX_s|^2 \rho_{\ve_{\ent{t}}}(\xi-\whX_s)\ud \xi \ud s 
    \\
    &+ 
    \frac{2}t \int_0^t\int_{\{|\xi|>R_{\ent{s}}\}} |\whX_s|^2 \rho_{\ve_{\ent{t}}}(\xi-\whX_s)\ud \xi \ud s 
    \\
    \le&   2 \ve_{\ent{t}}^2 \esp{|\zeta|^2}+ 
     \frac{2}t \int_0^t\int_{\{|\xi|>R_{\ent{s}}\}} |\whX_s|^2 \rho_{\ve_{\ent{t}}}(\xi-\whX_s)\ud \xi \ud s 
\end{align*}
using the same arguments as in the previous step for the first integral.
For the second term, we compute 
\begin{align*}
    \frac{1}t \int_0^t\int_{\{|\xi|>R_{\ent{s}}\}} |\whX_s|^2& \rho_{\ve_{\ent{t}}}(\xi-\whX_s)\ud \xi \ud s \\\ &\le \frac{1}t \int_0^t\int_{\{|\xi|>R_{\ent{s}}\}} |\whX_s|^2 \rho_{\ve_{\ent{t}}}(\xi-\whX_s)\ud \xi \ud s 
    \\
    &\le \frac{1}t \int_0^t\frac{1}{(R_{\ent{s}}-\ve_{\ent{t}}\,\sqrt{d})^a} \int_{\{|\xi|>R_{\ent{s}}\}} |\whX_s|^{2+a}\rho_{\ve_{\ent{t}}}(\xi-\whX_s)\ud \xi \ud s
    \\
    &\le \frac{1}t \int_0^t\frac{2^a}{R_{\ent{s}}^a} |\whX_s|^{2+a} \ud s
\end{align*}
thanks to the same arguments used in the previous step.  As $R_n=2 \sqrt{d}n^{\kappa}$, we deduce
\begin{align}
    \frac{1}t \int_2^t\int_{\{|\xi|>R_{\ent{t}}\}} |\whX_s|^2 \rho_{\ve_{\ent{t}}}(\xi-\whX_s)\ud \xi \ud s \le \frac{C}t \int_2^ts^{1-\kappa a }(\log s)^{1+\eta} \frac{|\whX_s|^{2+a}}{s(\log s)^{1+\eta}} \ud s
\end{align}
for  $\eta>0$. Then we compute 
\begin{align*}
    \esp{\int_2^t\frac{|\whX_s|^{2+a}}{s (\log s)^{1+\eta}} \ud s}= \int |y|^{2+a}\hm(\ud y)\int_2^t\frac{1}{s (\log (s))^{1+\eta}} \ud s < +\infty
\end{align*}
so that $\P(\int_2^t\frac{|\whX_s|^{2+a}}{s (\log (s))^{1+\eta}} \ud s< +\infty) =1$. Invoking the Kronecker Lemma~\ref{le kronecker gil}$(a)$, we obtain 
\begin{align}
    \frac{1}t \int_0^t\int_{\{|\xi|>R_{\ent{t}}\}} |\whX_s|^2 \rho_{\ve}(\xi-\whX_s)\ud \xi \ud s = t^{-\kappa a}o\big((\log t)^{1+\eta}\big) 
\end{align}
$\blacktriangleright$  We conclude that 
\begin{align*}
%\label{eq ccl I21}
    I_{21}(t) \le C\ve_{\ent{t}}^2+ t^{-\kappa a}o\big((\log t)^{1+\eta}\big) 
\end{align*} 
which, combined with~\eqref{eq I22} yields for $\eta>0$ 
\begin{align}\label{eq ccl I2}
    I_{2}(t) \le Ct^{-2r}+ t^{-\kappa a}o\big(\log^{1+\eta}(t)\big).
\end{align} 

\noindent {\sc Step~4} ({\em Conclusion}). Plugging the above bounds for $I_i(R)$, $i=1,2$, into~\eqref{eq:decomp1 as unbounded} and then in~\eqref{eq:decomp0 as unbounded} yields  for every $\eta>0$  small enough 
\begin{align*}
    {\cal W}_2(\hn_t, \hm)^2 
    %&\le  \textcolor{green}{C\big(\ve^2 + R^{-a} + R^{d+2} t^{-1/2}\ve^{-(d+1)}\big)}  \\
    % &
    % \le 
    % C\left( \frac{R_{\ent{t}}^{d+2}}{t\ve_{\ent{t}}^{d}} + \frac{R_{\ent{t}}^{d+2}}{\ve_{\ent{t}}^{d+1}} \ent{t}^{-\frac12} o\big( \log^{\frac 12+\eta} \ent{t}\big) + C\ve_{\ent{t}}^2+ t^{-\kappa a}o\big(\log^{1+\eta}(t)\big) +  6\ve_{\ent{t}}^2 \|\zeta\|^2_2 \right)
    % \\
     \le 
    C &\Big( t^{(d+3)\kappa + \textcolor{myred}{r(d+1)}-p} \\&+ 
    t^{(d+2)\kappa + \textcolor{jfc}{rd}+\theta-1}o(\log(t)^{\frac1{a+2} + \eta}) 
    +  t^{-2r} +  t^{-\kappa a}o\big(\log^{1+\eta}(t)\big)  \Big)
\end{align*}
for some real constant $\textcolor{black}{C=  C_{B,\Sigma,\rho,a,\eta}}\!\in (0,+\infty)$. 
Now, we seek to balance the various error terms. This amounts first to set $r=\kappa \frac{a}2$ and then to solve the linear system 
%\begin{align*}
%    -\kappa a &= (d+3)\kappa + r(d+1)-p 
%\\
% -\kappa a &= (d+2)\kappa + r(d+1)+\theta-1 
%\end{align*}
%which leads to
\begin{align*}
    \begin{cases}
    -\kappa a = (d+3)\kappa + \textcolor{myred}{\kappa \frac{a}2(d+1)}-p 
    \\
     -\kappa a = (d+2)\kappa + \textcolor{jfc}{\kappa \frac{ad}2} -\frac12  + \frac{1+pd}{a+2}
    \end{cases}
    \end{align*}
(having in mind that  $\theta-1 = -\frac12  + \frac{1+pd}{a+2}$) whose solution is given by 
% \begin{align*}
%     p  = \kappa (d+3)(1+a/2) \;\text{ and } \;\kappa = \frac{a}{2(a+2) \{2(d+2) + (d+3)(a+d) \}}.
% \end{align*}
\color{myred}
\begin{align*}
\kappa
=
\frac{a}
{(a+2)\left(d(d+3)+(a+2)(d+2)\right)}\;\text{ and }\; p
=
\frac{a(d+3)}
{2\left(d(d+3)+(a+2)(d+2)\right)}
\end{align*}

We observe that
\begin{align*}
\frac{2p}{a}=
\frac{d+3}
{d(d+3)+(a+2)(d+2)}<\frac1d
.
\end{align*} 
Eventually, we get 
\begin{align*}
    \cW_2(\hn_t,\hm) = o_\eta\big(t^{-\frac{a^2}
{2(a+2)\left(d(d+3)+(a+2)(d+2)\right)}}\log^{\frac12+\eta}(t)\big).
\end{align*}
%\color{black}
%
% \begin{align*}
%     \frac{2p}{a}=\frac{(d+3)}{ \{2(d+2) + (d+3)(a+d) \}}<\frac1d.
% \end{align*}
% Eventually, we get 
% \begin{align*}
%     \cW_2(\hn_t,\hm) = o_\eta\big(t^{-\frac{a^2}{2(a+2) \{2(d+2) + (d+3)(a+d) \}}}\log^{\frac12+\eta}(t)\big).
% \end{align*}
% \color{magenta}
% Now 
% \[
% \kappa = {\frac{1}{2(d+2)+a(d+3)}}\quad \mbox{ and }\quad r = \frac{\kappa a}{2}={\frac{a}{2(2(d+2)+a(d+3))}}.
% \]
% \textcolor{red}{We  check that $r(d+1) < \frac12$ so that $\theta >0$} and elementary computations show that 
% \[
% {\cal W}_2(\hn_t, \hm)\big) =o\big(t^{-\frac{a}{2((d+2)+a(d+3))}}\log^{1+\eta}(t)\big).
% \]
% for every small enough $\eta > 0$.
% Note that, as $a\to +\infty$, $\frac{a}{2((d+2)+a(d+3))}\to \frac{1}{2(d+3)}$.
% %\textcolor{magenta}{le meilleur taux avec cette preuve est $\frac{1}{d+1}$ à la limite}
% %}
\eproof
\color{black}

\bigskip

\color{black}
\noindent In the previous result, we observe that \textcolor{myred}{$\frac{a^2}
{2(a+2)\left(d(d+3)+(a+2)(d+2)\right)}\to \frac{1}{2(d+2)}$} 
%${\frac{a^2}{2(a+2) \{2(d+2) + (d+3)(a+d) \}}}\to \frac{1}{2(d+3)}$  
as $a\to +\infty$. This motivates our next proposition which is devoted to the case when the stationary distribution has an exponential moment. 
\color{black}

\begin{Proposition}\label{pr W2 as conv exp moment}
    Assume that $\Sigma$ is bounded and uniformly elliptic then, for $\eta>0$,
    \begin{align}\label{eq pr conv as exp moment}
        \cW_{2}(\hn_t,\hm)  = t^{-\textcolor{jfc}{\frac{1}{2(d+2)}}}o_\eta\left(\log(t)^{\frac12+\eta}\right).
        %O \left(
        %(\log t)^{\frac{d+2}{d+3}}t^{-\frac1{2(d+3)}} \right).
    \end{align}
    %\textcolor{red}{vérifier le terme en log}
\end{Proposition}

\color{black}
\proof 
Firstly, as $\Sigma$ is bounded, we know from Lemma~\ref{le sigma bounded exp moments}, that 
\begin{align} \label{eq ass exp moment}
    \exists\, \lambda>0\; \mbox{ such that }\; \int e^{\lambda |y|^2}  \hm(dy)<+\infty.  
\end{align}
We  will follow the steps,  lines and notations of the proof of Proposition~\ref{pr W2 as conv}.  
%\textcolor{red}{preuve à factoriser...}
%\\

\noindent {\sc Step~1} ({\em Decomposition}). 
% Let $\rho $ be a positive mollifier with support   $[-1,1]^d$ and let $\zeta$ be a   $\rho(y)\lambda_d(dy)$-distributed random vector. Let $\ve>0$ and let  us define   
%  $$
%  \rho_{\ve}(u) = \ve^{-d}\rho\Big(\frac{u}{\ve}\Big).
%  $$
% Note that $[\rho_{\ve}]_{\rm Lip} = \ve^{-(d+1)}[\rho]_{\rm Lip}=  \ve^{-(d+1)}\|\nabla \rho\|_{\sup}$ and $\|\rho_{\ve} \|_{\sup} \le \ve^{-d}\|\rho \|_{\sup}$. 
%
As in Step~1   of the proof of Proposition~\ref{pr W2 as conv},  we introduce a sequence of smoothing parameters $(\varepsilon_n)_{n \ge 1}$, here directly parametrized by 
\begin{align}\label{eq de epsilon}
    \textcolor{jfc}{
    \varepsilon_n = n^{-\frac{1}{2(d+2)}}\big(\tfrac{\log n}{d+2}\big)^{\frac12}.}
    % \\
    % \textcolor{yellow}{\epsilon_n = n^{-\frac{1}{2(d+3)}}(\log(n^\frac1{d+3}))^{\frac{d+2}{2(d+3)}}.}
\end{align}
%
% To control ${\cal W}_2(\hn_t, \hm)$ we  still rely on~\eqref{eq decomp Wasserstein unbounded} obtained by   the triangle inequality
% %\begin{align}\label{eq decomp Wasserstein}
% %{\cal W}_2(\hn_t, \hm) &\le {\cal W}_2(\hn_t, \hn_t*[\ve_{\ent{t}} \zeta]) +  {\cal W}_2(\hm, \hm *[\ve_{\ent{t}}\zeta])+{\cal W}_2(\hn_t *[\ve_{\ent{t}}\zeta], \hm *[\ve_{\ent{t}}\zeta]) .
% %\end{align}
% %From now on, we work assuming that $t \ge 1$.\\
% involving 
% %\noindent One checks that  the distributions $\hn_t*[\ve_{\ent{t}} \zeta]$ and $ \hm *[\ve_{\ent{t}}\zeta]$ have densities respectively  given  by 
% \[
% \frac{d\hn_t *[\ve_{\ent{t}}\zeta]}{d\lambda_d}(\xi) =   \frac{1}{t}\int_0^t\rho_{\ve_{\ent{t}}}(\xi-\whX_s) ds
% \quad\mbox{ and }\quad 
% \frac{d\hm_t *[\ve_{\ent{t}}\zeta]}{d\lambda_d}(\xi) = \int_{\R^d}\rho_{\ve_{\ent{t}}}(\xi-y)\hm(dy).  
% \]
% and set 
% \begin{equation*}
% %\label{eq:dfnphieps as}
% \phi_{\ve_{\ent{t}}}(\xi, x) =\rho_{\ve_{\ent{t}}}(\xi-x)-\int_{\R^d}\rho_{\ve_{\ent{t}}}(\xi-y)\hm(dy) 
% \end{equation*}
% then, it follows from~\eqref{eq decomp Wasserstein unbounded} that  
% \begin{align} \label{eq:decomp0 as}
%  {\cal W}_2(\hn_t, \hm)^2
% &\le 6\ve_{\ent{t}}^2 \|\zeta\|^2_2 +3 \int_{\R^d} |\xi|^2  \Big|\frac1t\int_0^t \phi_{\ve_{\ent{t}}} (\xi, \whX_s) ds \Big| \ud \xi.   
% %& \le 6\ve^2 \|\zeta\|^2_2 +3\bigg(\int_{\R^d} |\xi|^2  \E\, \Big|\frac1t\int_0^t \phi_{\ve} (\xi, X_s) ds \Big| \ud \xi \bigg),
% \end{align}
% %where we used $\|\cdot\|_1\le \|\cdot\|_2$ and Fubini's theorem in the second line.
%
Invoking Lemma \ref{le smoothing generic} (with the above smoothing parameter sequence), we have from \eqref{eq:decomp0 le} that
\begin{align}\label{eq:decomp0 as unbounded}
 {\cal W}_2(\hn_t, \hm)^2
&\le 6\ve_{\ent{t}}^2 \|\zeta\|^2_2 + 9 \int_{\R^d} |\xi|^2  \Big|\frac1t\int_0^t \phi_{\ve_{\ent{t}}} (\xi, \whX_s) ds \Big| \ud \xi.   
\end{align}

We consider here the following increasing sequence $(R_n)_{n \ge 0}$ of truncation parameters   
\begin{align}\label{eq de R}
    R_n &= 2\Big(\frac{\log n}{(\lambda \wedge 1)\textcolor{jfc}{(d+2)}}\Big)^{\frac12} , \quad n\ge 2, \text{ and } R_0:=R_1 := R_2.
\end{align}
  In the previous definition, $\lambda$ is given by  Lemma~\ref{le sigma bounded exp moments}. 
  We then have that, for  some positive integer $N_{d}\ge 2$,
\begin{align}\label{eq small truc pour plus tard bis}
    R_n \ge 2 \sqrt{d}\,\varepsilon_{n'}, \;\text{ for } \;n,n'\ge N_{d}.
\end{align}
This will prove useful at some point in the proof below. 

%  \textcolor{magenta}{\bf verifier pourquoi on part de $2\epsilon \sqrt{d}$ pour la dernière partie de la preuve...}:
% \begin{align*}
%     \R_+ \ni t \mapsto R(t) \in (2\epsilon \sqrt{d},+\infty),
% \end{align*}
% such that $R(t)\sim t^\kappa$, for some $\kappa > 0$ to be chosen later on.
% \textcolor{red}{bien reprendre parametrisation...}
\noindent Now, we inspect  the terms $I_1(t)$, $i=1,2$ as defined in the poof of Proposition~\ref{pr conv L2} according to \eqref{eq:decomp1 as unbounded}.
%split the  above integral into two terms using $R$ as follows: 
%\begin{equation}\label{eq:decomp1 as}
% \int_{\R^d} |\xi|^2\E\, \Big|\frac1t\int_0^t \phi_{\ve(t)} (\xi, \whX_s) ds \Big| \ud \xi  = I_1(t)  +I_2(t) 
%\end{equation} 
% with
% \[
%I_1(t)=  \int_{\{|\xi| \le R_{\ent{t}}\}} |\xi|^2\Big|\frac1t\int_0^t \phi_{\ve_{\ent{t}}} (\xi, \whX_s) ds \Big| \ud \xi\; \mbox{ and }\; I_2(t)= \int_{\{|\xi|> R_{\ent{t}}\}} |\xi|^2 \Big|\frac1t\int_0^t \phi_{\ve_{\ent{t}}} (\xi, \whX_s) ds \Big| \ud \xi .
%\]

\smallskip
\noindent {\sc Step 2} ({\em Term $I_1(t)$}). 
%\textcolor{yellow}{I do the computation with $d=1$ to start with.} 
% Let $(\gamma_n)$ a decreasing sequence of steps and $\Gamma_n = \sum_{k=1}%^n\gamma_k$, it is assumed that $\Gamma_n \rightarrow +\infty$ when $n \rightarrow %+\infty$. Recall notation: $\t$, $N(t)$.
% \textcolor{red}{on peut prendre un pas de temps constant !! $\Gamma_n \sim n$ %simplement donc...}
%
We know  that 
\begin{align*}
    I_1(t)
    %&
    %\le \int_{|\xi|\le R_{\ent{t}}}|\xi|^2\Big |\frac1{t}\int_{\ent{t}}^t\phi_{\ve_{\ent{t}}}(\xi,\whX_s)\ud s\Big|\ud \xi + \int_{|\xi|\le R_{\ent{t}}}|\xi|^2 \Big|\frac1{\ent{t}}\int_{0}^{\ent{t}}\phi_{\ve_{\ent{t}}}(\xi,\whX_s)\ud s\Big|\ud \xi
    %\\
    &\le C\frac{R_{\ent{t}}^{d+2}}{t\ve_{\ent{t}}^{d}}+ \tilde{I}_1(\ent{t})
    % \int_{|\xi|\le R(\Gamma_{N(t)+1})}|\xi|^2 |\frac1{\Gamma_{N(t)}}\int_{0}^{\Gamma_{N(t)}}\phi_{\ve(t)}(\xi,X_s)\ud s|\ud \xi
    % =:I_{12}(t)
\end{align*} 
where
%, for $n \ge 1$, 
\begin{align*}
    \tilde{I}_1(n) &:= \int_{|\xi|\le R_n}|\xi|^2 |\frac1{n}\int_{0}^{n}\phi_{\ve_n}(\xi,\whX_s)\ud s|\ud \xi\\
%
%%
%%We divide $[-R_n,R_n]$ into $K_n$ interval of size $\delta_n := 2R_n/K_n$ and subsequently $[-R_n,R_n]^d$ is covered with $K_n^d$ hypercube of length $\delta_n$. For a given indexation $k=0,\dots,K_n^d-1$, we denote these hypercubes $A^n_k$.
%%Then,
%\begin{align*}
%    \tilde{I}_{1}(n) \le R_n^2 \sum_{k=0}^{K_{n}^d-1}\int_{\xi \in A^n_k} \Big|\frac1{n}\int_{0}^{n}\phi_{\ve_n}(\xi,\whX_s)\ud s\Big|\ud \xi
%\end{align*}
%For $\xi \in A^n_k$, we observe that
%\begin{align}
%    |\phi_{\ve_n}(\xi,\whX_s)-\phi_{\ve_n}(\xi^n_k,\whX_s)| \le C\frac{\delta_n}{\ve_n^{d+1}}
%\end{align}
%where $\xi^n_k$ is a point in $A^n_k$ (the lower corner of $A^n_k$ say).
%This leads to, observing that $\int_{\xi \in A^n_k}\ud \xi \le C\delta_n^d$,
%\begin{align}
%  &\le C \frac{R_n^{d+2}}{\ve_n^{d+1}} \delta_n + R_n^2 \delta_n^d \sum_{k=0}^{K_{n}^d-1} |\frac1{n}\int_{0}^{n}\phi_{\ve_n}(\xi^n_k,\whX_s)\ud s|
%\end{align}
%and 
% \begin{align}
 % \tilde{I}_{1}(n) 
   &\le C \frac{R_n^{d+2}}{ \textcolor{myred}{\ve_n^{d+1}}} \delta_n + R_n^{d+2} \sup_{k=0,\ldots,K_n^d-1}\Big|\frac1{n}\int_{0}^{n}\phi_{\ve_n}(\xi^n_k,\whX_s)\ud s\Big|.
 \end{align*}
Now we set $\psi(n):=\sqrt{n\log(1+n)}$ and following the lines of the proof of Proposition~\ref{pr W2 as conv}, we get 
%Now, invoking Proposition~\ref{le unif elliptic cobord}, we have 
%\begin{align*}
%   \Big|\frac1{n}\int_{0}^{\n}\phi_{\ve_n}(\xi^n_k,\whX_s)\ud s\Big|
%%    &\le \frac1{n}|\omega_{\ve_n}(\xi^n_k,\whX_{n})-\omega_{\ve_n}(\xi^n_k,\whX_0)|
%%    + \frac1{n}|\int_{0}^{n}\partial_x \omega_{\ve_n}(\xi^n_k,\whX_s)\Sigma(\whX_s)\ud W_s|
%%    \\
%    &\le  C\frac{|\whX_{n}-\whX_{0}|}{\ve_n^{d+1}n}
%    + \frac{\psi(n)}{n}\Big|\frac1{\psi(n)}\int_{0}^{n}\partial_x \omega_{\ve_n}(\xi^n_k,\whX_s)\Sigma(\whX_s)\ud W_s\Big|
%\end{align*}
%where $\psi(n):=\sqrt{n\log(1+n)}$.
%\\
%So that 
\begin{align*}
%\label{eq temp Itilde un}
    \tilde{I}_{1}(n) \le C \frac{R_n^{d+2}}{\textcolor{jfc}{\ve_n^{d}}} \left(\textcolor{myred}{\frac{\delta_n}{\ve_n}}+\frac{|\whX_{n}-\whX_{0}|}{n} +  \frac{\psi(n)}{{n}}\sup_{k=0,\ldots,K_n^d-1}\Big|\frac1{\psi(n)}\int_{0}^{n}\textcolor{jfc}{\ve_n^{d}}\partial_x \omega_{\ve_n}(\xi^n_k,\whX_s)\Sigma(\whX_s)\ud W_s\Big|\right)
\end{align*}

% \textcolor{red}{Je me suis arrêté là pour cette étape. En effet, quand on passe de~\eqref{eq temp Itilde un} à~\eqref{eq majo violente} il y a $\frac1{\ve_n^{d+1}}$ qui est mis en factor et donc pour le terme martingale ci-dessus on devrait avoir $\ve_n^{d+1}$ qui compense $|\partial_x \omega_{\ve_n}|_\infty$. Du coup j'ai eu l'impression qu'on avbait un peu plus de place pour obtenir une vitesse de convergence (toute pourrie) dans le cas sigma non-bornée... je me suis endormi en essayant d'optimiser les paramètres... je n'ai donc pas fini de vérifier l'étape 2 et donc je n'ai pas repris l'étape 4 de conclusion.}
\noindent Now we study the last term in the r.h.s. of the previous inequality.
% \begin{align}\label{eq as conv gyorfi}
%     \sup_{k=0,\ldots,K_n^d-1}\Big|\frac1{\psi(n)}\int_{0}^{{n}}\ve^{d+1}_n\partial_x \omega_{\ve_n}(\xi^n_k,\whX_s)\Sigma(\whX_s)\ud W_s\Big|\rightarrow 0 \text{ as } n \rightarrow +\infty.
% \end{align}
Still using 
\[
Z^{n,k}_s =    \textcolor{jfc}{\ve_n^{d}}\partial_x \omega_{\ve_n}(\xi^n_k,\whX_s)\Sigma(\whX_s)
\]
we observe that, as $\Sigma$ is bounded and $\omega_{\epsilon_n}$ is Lipschitz,
\begin{align}\label{eq control integral}
    |Z^{n,k}_s|  \le \mathfrak{c}. 
\end{align}
Let $\ell>\mathfrak{c}\sqrt{2(pd+1)}$ be a rational number. Applying~\eqref{eq thanks hoeff} in Lemma~\ref{le hoeffding like}$(a)$ yields 
\begin{align*}
    \P\left(\Big|\frac1{\psi(n)}\int_{0}^{n}Z^{n,k}_s\ud W_s\Big| > \ell\right)
    \le 2e^{ -\frac{\ell^2}{2\mathfrak{c}^2}n^{-1}\psi(n)^2}.
\end{align*}

Then,
\begin{align*}
    \P\left(\sup_{k=0,\ldots,K_n^d-1}\Big|\frac1{\psi(n)}\int_{0}^{{n}}Z^{n,k}_s\ud W_s\Big| > \ell\right)
    \le 2K_n^d e^{ -\frac{\ell^2}{2\mathfrak{c}^2}n^{-1}\psi(n)^2}.
\end{align*}
We now set $K_n = \lfloor n^{p}\rfloor\sim n^{p}$ for some $p$ to be determined later on. From the definition of $\psi$, we  compute
\begin{align*}\
%label{eq majo sans scrupule bounded}
    \sum_{n=1}^{+\infty}\P\left(\sup_{k=0,\ldots,K_n^d-1}\Big|\frac1{\psi(n)}\int_{0}^{n}Z^{n,k}_s\ud W_s\Big| > \ell\right)
    \le 2 \sum_{n=1}^{+\infty}   e^{ (pd-\frac{\ell^2}{2\mathfrak{c}^2}) \log n}
    = 2 \sum_{n=1}^{+\infty} n^{(pd-\frac{\ell^2}{2\mathfrak{c}^2})}.
    % < +\infty
\end{align*}
It follows that the series in the righthand side of the above equation is converging given our choice for $\ell$, so that, by the Borel-Cantelli Lemma, 
\begin{align*}
%\label{eq you deserve it0}
    \limsup_n \sup_{k=0,\ldots,K_n^d-1}\Big|\frac1{\sqrt{n\log n}}\int_{0}^{{n}}Z^{n,k}_s \ud W_s\Big| \le \ell \quad a.s.
\end{align*}
or, equivalently
\begin{align*}
%\label{eq you deserve it}
    \limsup_n \sup_{k=0,\ldots,K_n^d-1} \Big|\frac{1}{\sqrt{n\log n}}\int_{0}^{{n}}Z^{n,k}_s\ud W_s\Big| \le \mathfrak{c}\sqrt{2(pd+1)} \quad a.s.
\end{align*}

Furthermore, it follows from Lemma~\ref{le control SDE}$(b)$ (with $a=2$)  that
\begin{align*}
    \frac{|X_{n}-X_{0}|}{n} = o_\eta\big(n^{-\frac34}(\log n)^{\frac14+\eta}\big) = o\Big( \frac{\psi(n)}{n}\Big).
\end{align*}
%Consequently this term can be neglected in regards to $\frac{\psi(n)}{n}$.

This leads to
\begin{align*}
%\label{eq majo i tilde un bounded}
    \tilde{I}_{1}(n) &\le C \frac{R_n^{d+2}}{\textcolor{jfc}{\ve_n^{d}}} \left(\textcolor{myred}{\frac{\delta_n}{\ve_n}}
    + o_\eta\big(n^{-\frac12}(\log n)^{\frac12+\eta}\big) 
    \right) .
    % \\
    % & \le C \frac{R_n^{d+2}}{\ve_n^{d+1}} \left(\delta_n
    % + o\big(n^{-\frac12}(\log n)^{\frac12+\eta}\big) \right)
\end{align*}

Since $R_n \sim c \sqrt{\log n}$ and $\delta_n=\frac{2R_n}{K_n}$ we can set \textcolor{myred}{$p = \frac12-\frac1{2(d+2)}$}, to obtain finally
%
%we obtain 
\begin{align*}
%\label{eq majo almost done bounded sig as}
    \tilde{I}_{1}(n) = \frac{R_n^{d+2}}{\textcolor{jfc}{\ve_n^{d}}} 
     o_\eta\big(n^{-\frac12}(\log n)^{\frac12+\eta}\big)\quad\mbox{ for every $\eta > 0$}. 
\end{align*}

\smallskip
\noindent {\sc Step~3} ({\em Term $I_2$}). We will successively inspect the two terms $I_{21}(t)$ and $I_{22}(t)$ of the decomposition~\eqref{eq:decompI_2} of $I_2(t)$.
%We remark that
%%, using Fubini's Theorem at each line  and the elementary inequality $(a+b)^2 \le 2(a^2+b^2)$, $a,b\ge 0$, in the third line that 
%\begin{align*}
%% I_2(R(t)) &\le \int_{\{|\xi|>R(t)\}}|\xi|^2\bigg( \frac 1t \int_0^t\rho_{\ve}(\xi-X_s)ds + \int \rho_{\ve}(\xi-y)\hm(dy)\bigg) \ud \xi\\
%I_2(t) &\le 
%\int_{\{|\xi|>R_{\ent{t}}\}}\frac{|\xi|^2}t \int_0^t\rho_{\ve_{\ent{t}}}(\xi-\whX_s)\ud s \ud \xi 
%+ \int_{\{|\xi|>R_{\ent{t}}\}}|\xi|^2\int \rho_{\ve_{\ent{t}}}(\xi-y)\hm(\d y)\ud \xi
%=:I_{21}+I_{22} \\
%% &\le 2\int_{\{|\xi|>R(t)\}} |\xi|^2 \bigg(\frac1t \int_0^t\rho_{\ve}(\xi-X_s)ds\bigg)\ud \xi\\
%% & \le  4\bigg[\frac1t\int_0^t \int_{\{|\xi|>R\}} \E\, |\xi-X_s|^2  \rho_{\ve}(\xi-X_s)\ud \xi ds + \frac1t \int_0^t\int_{\{|\xi|>R\}} \E\, |X_s|^2\rho_{\ve}(\xi-X_s)\ud \xi ds  \bigg].
%\end{align*}
%\textcolor{red}{les points 3.a et 3.b ci-dessous ont été mis à jour}\\

\noindent $\blacktriangleright$   We first study $I_{22}$ which itself can be upper-bounded as follows
\begin{align*}
    \int_{\{|\xi|>R_{\ent{t}}\}}|\xi|^2\int \rho_{\ve_{\ent{t}}}(\xi-y)\hm(\ud y)\ud \xi
    \le& 2   \int_{\{|\xi|>R_{\ent{t}}\}}\int |\xi-y|^2\rho_{\ve_{\ent{t}}}(\xi-y)\hm(\ud y)\ud \xi 
    \\
    &+ 2 \int_{\{|\xi|>R_{\ent{t}}\}}\int |y|^2\rho_{\ve_{\ent{t}}}(\xi-y)\hm(\ud y)\ud \xi.
\end{align*}
The first term in the r.h.s. still satisfies 
\begin{align}\label{eq:>Rdustep2}
    \int_{\{|\xi|>R_{\ent{t}}\}}\int |\xi-y|^2\rho_{\ve_{\ent{t}}}(\xi-y)\hm(\ud y)\ud \xi 
%    &\le
%    \int\int  |\xi-y|^2\rho_{\ve_{\ent{t}}}(\xi-y)\hm(\ud y)\ud \xi
%    \\
    &
    = |z|^2\rho_{\ve_{\ent{t}}}(z)\ud z = \ve_{\ent{t}}^2 \esp{|\zeta|^2}
\end{align}
and  the second term is also still upper-bounded by 
\begin{align*}
    \int_{\{|\xi|>R_{\ent{t}}\}}\int |y|^2 \rho_{\ve_{\ent{t}}}(\xi-y)\hm(\ud y)\ud \xi
    \le   
    \int_{\{|\xi|>R_{\ent{t}}\}}\int |y|^2\1_{\set{|y-\xi|<\sqrt{\ve_{\ent{t}}}d}}\rho_{\ve_{\ent{t}}}(\xi-y)\hm(\ud y)\ud \xi.
\end{align*}
Still using  that 
$$
\set{y:| y-\xi|\le \ve_{\ent{t}}\sqrt{d}, |\xi|> R_{\ent{t}}}\subset \set{y:|y|\ge R_{\ent{t}}- \ve_{\ent{t}}\sqrt{d}, |\xi|> R_{\ent{t}}}.
$$
For $t \ge N_d$, $R_{\ent{t}} \ge 2\sqrt{d}\ve_{\ent{t}}$, according to~\eqref{eq small truc pour plus tard bis}, we then have 
\begin{align*}
\set{|y-\xi|\le \ve_{\ent{t}}\sqrt{d}, |\xi|> R_{\ent{t}}}\subset 
   \set{|y|\ge \frac{R_{\ent{t}}}2, |\xi|> R_{\ent{t}}}.
\end{align*}
It follows then by combining exponential Markov inequality with  Lemma~\ref{le sigma bounded exp moments}  that 
% \color{yellow}
% \begin{align}
%     I_{22}(R) & \le 
%     \esp{e^{\frac{\lambda}2|\whX_s|^2}|\whX_s|^2e^{-\frac{\lambda R^2}4}\int_{\{|\xi|>R\}}\rho_{\ve}(\xi-\whX_s) \ud \xi }
%     \\
%     &\le C_\lambda e^{-\frac{\lambda R^2}4} \esp{e^{\lambda|\whX_s|^2}}
%     \\
%     &\le C_\lambda e^{-\frac{\lambda R^2}4}.
% \end{align}
% \color{blue}
\begin{align*}
 \int_{\{|\xi|>R_{\ent{t}}\}}\int \mbox{\bf1}_{\{|y-\xi|\le\sqrt{d}\, \ve_{\ent{t}}\}} &|y|^2\rho_{\ve_{\ent{t}}}(\xi-y)\hm(\ud y)\ud \xi\\
 &\le 
 \int_{\{|\xi|>R_{\ent{t}}\}} e^{-\frac{\lambda R_{\ent{t}}^2}4} \int|y|^{2}e^{\frac{\lambda}2 |y|^2}\rho_{\ve}(\xi-y)\hm(\ud y)\ud \xi\\
 &\le 
 e^{-\frac{\lambda R_{\ent{t}}^2}4}   \int|y|^{2}e^{\frac{\lambda}2 |y|^2}\int\rho_{\ve}(\xi-y) \ud \xi \hm(\ud y) \\
 &\le 
 e^{-\frac{\lambda R_{\ent{t}}^2}4}   \int|y|^{2}e^{\frac{\lambda}2 |y|^2} \hm(\ud y) \\
 &\le 
 e^{-\frac{\lambda R_{\ent{t}}^2}4}  C_\lambda \int e^{\lambda |y|^2} \hm(\ud y) 
 \\
    &\le C_\lambda e^{-\frac{\lambda R_{\ent{t}}^2}4}
 %&\le \frac{1}{(R_{\ent{t}}-\ve_{\ent{t}}\,\sqrt{d})^a}  \int|y|^{2+a}\int \rho_{\ve_{\ent{t}}}(\xi-y)\ud \xi\mu^\star(\ud y)\\
%&\le \frac{1}{(R_{\ent{t}}-\ve\,\sqrt{d})^a}  \int|y|^{2+a}\mu^\star(\ud y)
%\\
 %&\le \frac{2^a}{R_{\ent{t}}^a}\int |y|^{2+a}\hm(\ud y)
\end{align*}
so that, finally, 
%according to our choice of $(R_n)-_n$ we have
\begin{align}\label{eq I22 bis}
    I_{22}(t) \le C\big(\ve^2_{\ent{t}} + e^{-\frac{\lambda R_{\ent{t}}^2}4}\big).
\end{align}
\noindent $\blacktriangleright$   We now study $I_{21}$.
We first observe, since $R_n$ is non-decreasing, that
\begin{align*}
    \int_{\{|\xi|>R_{\ent{t}}\}}\frac{|\xi|^2}t \int_0^t\rho_{\ve_{\ent{t}}}(\xi-\whX_s)\ud s \ud \xi
%    \le &
%    \frac{2}t \int_0^t\int_{\{|\xi|>R_{\ent{t}}\}} |\xi - \whX_s|^2 \rho_{\ve_{\ent{t}}}(\xi-\whX_s)\ud \xi \ud s 
%    \\
%    &+ 
%    \frac{2}t \int_0^t\int_{\{|\xi|>R_{\ent{s}}\}} |\whX_s|^2 \rho_{\ve_{\ent{t}}}(\xi-\whX_s)\ud \xi \ud s 
%    \\
    &\le    2 \ve_{\ent{t}}^2 \esp{|\zeta|^2}+ 
     \frac{2}t \int_0^t\int_{\{|\xi|>R_{\ent{s}}\}} |\whX_s|^2 \rho_{\ve_{\ent{t}}}(\xi-\whX_s)\ud \xi \ud s 
\end{align*}
with the same arguments as in~\eqref{eq:>Rdustep2} for the first integral.
For the second one, one has
\begin{align}
    \nonumber \frac{1}t \int_0^t\int_{\{|\xi|>R_{\ent{s}}\}} |\whX_s|^2& \rho_{\ve_{\ent{t}}}(\xi-\whX_s)\ud \xi \ud s 
    \\
       \nonumber &\le 
    \frac{1}t \int_0^{N_d} |\whX_s|^2 \ud s 
    +
    \frac{1}t \int_{N_d}^t\int_{\{|\xi|>R_{\ent{s}}\}} |\whX_s|^2 \rho_{\ve_{\ent{t}}}(\xi-\whX_s)\ud \xi \ud s 
    \\
       \nonumber &\le O(t^{-1}) 
    + \frac{1}t \int_{N_d}^te^{-\frac{\lambda R^2_{\ent{s}}}4} \int_{\{|\xi|>R_{\ent{s}}\}} |\whX_s|^{2}e^{\frac{\lambda}2 |\whX_s|^2}\rho_{\ve_{\ent{t}}}(\xi-\whX_s)\ud \xi \ud s
    \\
       \nonumber &\le O(t^{-1}) 
    + \frac{1}t \int_{N_d}^te^{-\frac{\lambda R^2_{\ent{s}}}4}  |\whX_s|^{2}e^{\frac{\lambda}2 |\whX_s|^2} \ud s
    \\
    &\le O(t^{-1}) 
    + \frac{C_\lambda}t \int_{N_d}^te^{-\frac{\lambda R^2_{\ent{s}}}4} e^{{\lambda} |\whX_s|^2} \ud s,
\end{align}
where we used the same arguments as in Step~2.  According to our choice of $(R_n)$ in~\eqref{eq de R}, we have then 
\begin{align*}
    \frac{1}t \int_{N_d}^te^{-\frac{\lambda R^2_{\ent{s}}}4} e^{{\lambda} |\whX_s|^2} \ud s
    \le \frac{C}t \int_{N_d}^t s^{1-\textcolor{jfc}{\frac{1}{d+2}}} (\log s)^{1+\eta} \frac{e^{{\lambda} |\whX_s|^2}}{s(\log s)^{1+\eta}}\ud s
\end{align*} 
for $\eta>0$. Then we compute 
\begin{align}
    \esp{\int_{N_d}^{+\infty} \frac{e^{{\lambda} |\whX_s|^2}}{s(\log s)^{1+\eta}}\ud s}
 = \int e^{\lambda |y|^2}\hm(\ud y)\int_{N_d}^{+\infty} \frac{1}{s(\log s)^{1+\eta}}\ud s<+\infty
\end{align}
using Lemma~\ref{le sigma bounded exp moments}.
Thus $ \int_{N_d}^{+\infty} \frac{e^{{\lambda} |\whX_s|^2}}{s\log(s)^{1+\eta}}\ud s<+\infty \;\P$-$a.s.$ and   Kronecker's Lemma~\ref{le kronecker gil}$(a)$ implies
\begin{align}
    \frac{1}t \int_{N_d}^te^{-\frac{\lambda R^2_{\ent{s}}}4} e^{{\lambda} |\whX_s|^2} \ud s = t^{-\textcolor{jfc}{\frac{1}{d+2}} }o\big((\log t)^{1+\eta}\big).
\end{align}
\noindent $\blacktriangleright$  We conclude that 
\begin{align}\label{eq ccl I21 bis}
    I_{21}(t) \le C\ve_{\ent{t}}^2+ t^{-\textcolor{jfc}{\frac{1}{d+2}}  }o\big((\log t)^{1+\eta}\big) 
\end{align} 
which combined with~\eqref{eq I22 bis} yields
\begin{align}\label{eq ccl I2 bis}
    I_{2}(t) \le C\ve_{\ent{t}}^2+ t^{-\textcolor{jfc}{\frac{1}{d+2}}  }o\big((\log t)^{1+\eta}\big)\quad \mbox{for any $\eta>0$}.
\end{align}

\medskip
\noindent {\sc Step~4} ({\em Conclusion}). Collecting all the above bounds for $I_i(R)$, $i=1,2$  yields  for  $\eta>0$.
\begin{align*}
    {\cal W}_2(\hn_t, \hm)^2 
    %&\le  \textcolor{green}{C\big(\ve^2 + R^{-a} + R^{d+2} t^{-1/2}\ve^{-(d+1)}\big)}  \\
    &
    \le 
    C\left(  \frac{R_{\ent{t}}^{d+2}}{\textcolor{jfc}{\ve_{\ent{t}}^{d}}} \ent{t}^{-\frac12} o\big( \log^{\frac 12+\eta} \ent{t}\big) + \ve_{\ent{t}}^2+ t^{-\textcolor{jfc}{\frac{1}{d+2}}  }o\big(\log^{1+\eta}(t)\big)   \right)
    % \\
    % & \le 
    % C \left( \frac{R_{\ent{t}}^{d+2}}{\ve_{\ent{t}}^{d+1}} \ent{t}^{-\frac12} o\big( \log^{\frac 12+\eta} \ent{t}\big) + \ve_{\ent{t}}^2+  t^{-\kappa a}o\big(\log^{1+\eta}(t)\big)  \right)
    \\
    &= t^{-\textcolor{jfc}{\frac{1}{d+2}}}o\left(\log(t)^{1+\eta}\right).\hskip 7cm \Box
\end{align*}
% for some real constant $\textcolor{blue}{C=  C_{b,\s,\rho,a,\eta}}\!\in (0,+\infty)$. Now we can set 
% \[
% \kappa = {\frac{1}{2(d+2)+a(d+3)}}\quad \mbox{ and }\quad r = \frac{\kappa a}{2}={\frac{a}{2(2(d+2)+a(d+3))}}.
% \]
% \textcolor{red}{We  check that $r(d+1) < \frac12$ so that $\theta >0$} and elementary computations show that 
% \[
% {\cal W}_2(\hn_t, \hm)\big) =o\big(t^{-\frac{a}{2((d+2)+a(d+3))}}\log^{1+\eta}(t)\big).
% \]
% for every small enough $\eta > 0$.
% Note that, as $a\to +\infty$, $\frac{a}{2((d+2)+a(d+3))}\to \frac{1}{2(d+3)}$.
% %\textcolor{magenta}{le meilleur taux avec cette preuve est $\frac{1}{d+1}$ à la limite}
% %} 
\color{black}

%% file: appendix.tex
\subsection{Proof of integrability properties: Lemma \ref{le sigma bounded exp moments}}

1. Combining $\HYP{L}$ and $\HYP{C}$, one obtains that Hajek's criterion holds, namely the  existence of $K'\ge 0$ and $\a'\!\in (0,\a)$ such that 
\begin{align}\label{eq Hajek for classic SDE}
  \forall\, x\!\in \R^d,\qquad 2(B(x)|x) + \| \Sigma(x)\|^2_{_F} \le K'-\a'|x|^2. 
\end{align}
Then, it is well known that the SDE \eqref{eq solution sde hat} has a unique invariant distribution $\hm \in \cP_{2+a}(\R^d)$  for some $a>0$.
\\
2. We now give the proof of the bounded exponential moments when $\Sigma$ is bounded.

% \subsection{Bounded exponential moments: proof of Lemma~\ref{le sigma bounded exp moments}}\label{subsec expo moment}
\color{black}
%\begin{Proposition} Assume \ref{eq Hajek for classic SDE} with mean-reverting coefficient $\a>0$ is in force and $\Sigma$is bounded. Then for every $\la \!\!\in \big(0,\frac{\a}{2  \|\Sigma\|_{\sup}} \big)$, then\
%%tiny \textcolor{red}{@Gil: pas fini:faut int\'egrer par rapport \`a $\hat \mu$ en $K$-localisant}\normalsize
%\[ 
% \int_{\R^d} e^{\la |x|^2} \hat \mu(\ud x)<+\infty.
%\]
%\end{Proposition}
%\noindent {\em Proof.} 
We denote by ${\cal L}$ the infinitesimal generator of the diffusion~\eqref{eq solution sde hat}  and  temporarily  set $V(x)= |x|^2+1$, $x\!\in \R^d$. We assume the mean-reverting assumption condition~\eqref{eq Hajek for classic SDE}, which reads, using the infinitesimal generator,
\[
{\cal L}(V) \le K- \textcolor{myred}{\tilde\a} V
\]
 where $\textcolor{myred}{\tilde\a} >0$ and $K\ge 0$. We also assume that $\Sigma$ is a bounded function.  
 One has, for $\la >0$
 \begin{align*}
{\cal L}(e^{\la V})(x) &= \la e^{\la V(x)}\Big(2(x |B(x)) + {\rm Tr}(\Sigma\Sigma^{\top}(x) ) + 2\la |\Sigma(x)x|^2\Big)\\
& \le \la e^{\la V(x)}\Big(K-\textcolor{myred}{\tilde\a}  V(x)   + 2\la \|\Sigma\|_{\sup}^{\textcolor{myred}{2}}V(x)\Big).
 \end{align*}
Consequently, if $\la <\frac{\textcolor{myred}{\tilde\a} }{2  \|\Sigma\|^{\textcolor{myred}{2}}_{\sup}} $, $\textcolor{myred}{\bar\a} = \textcolor{myred}{\tilde\a}  -2\la \|\Sigma\|^{\textcolor{myred}{2}}_{\sup}>0$ and 
\[
{\cal L}(e^{\la V})\le  e^{\la V}\big(K-\textcolor{myred}{\bar\a}V \big).
\]
Now, for such  $\la$, 
\[
  e^{\la V}\big(K-\textcolor{myred}{\bar\a}V \big) +\textcolor{myred}{\bar\a} e^{\la V} = e^{\la V} \big(K +\textcolor{myred}{\bar\a}(1-V)\big).
\]
The function $K +\textcolor{myred}{\bar\a}(1-V)$ is continuous  and goes to $-\infty$ at infinity so that $\textcolor{myred}{\bar K} = \sup_{\R^d} e^{\la V}(K +\textcolor{myred}{\bar\a}(1-V) )<+\infty$. This in turn implies that
\[
e^{\la V} (K -\textcolor{myred}{\bar\a} V) \le \textcolor{myred}{\bar K} -\textcolor{myred}{\bar\a} e^{\la V}.
\]
Now  let $x_0\!\in \R^d$ and set for every $L>0$, $\tau_L := \inf\{t  \ge 0: |X^{x_0}_t|\ge L\} $.
Applying It\^o's formula between $0$ and $t\wedge \tau_L$ to $e^{\textcolor{myred}{\bar \a}  t }e^{\la V(X^{x_0}_t)}$ yields
\[
e^{\textcolor{myred}{\bar \a}(t\wedge \tau_L)+\la V(X^{x_0}_{t\wedge \tau_L})}\le  e^{ \la V(x_0)}+\int_0^{t}e^{\textcolor{myred}{\bar \a} s} \textcolor{myred}{\bar K}\ud s+2 \la \int_0^{t\wedge\tau_L}e^{\textcolor{myred}{\bar \a} s+\la V(X^{x_0}_s)}\big(X^{x_0}_{s}\,|\,\Sigma(X^{x_0}_{s\wedge \tau_L})\big)\ud W_s
\]
so that, taking expectation,    
\[
\E\, e^{\textcolor{myred}{\bar \a}(t\wedge \tau_L) +\la V(X^{x_0}_{t\wedge \tau_L})}  \le   e^{ \la V(x_0)}+ \textcolor{myred}{\bar K}\frac{e^{\textcolor{myred}{\bar \a}t}-1}{\textcolor{myred}{\bar \a}}.
%+2 \la \int_0^{t\wedge\tau_L}e^{\a' s\la +V(X^{x_0}_{s}}\big(X^{x_0}_{s}\,|\,\sigma(X^{x_0}_{s}\ud W_s\big).
\]
By Fatou's Lemma applied on the l.h.s.  when $L\uparrow +\infty$, we derive 
\[
e^{\textcolor{myred}{\bar \a} t}\E\, e^{\la V(X^{x_0}_{t})} \le e^{ \la V(x_0)}+  \textcolor{myred}{\bar K}\frac{e^{\textcolor{myred}{\bar \a} t}-1}{\textcolor{myred}{\bar \a}}
\]
i.e., for every $t\ge 0$, 
%\textcolor{red}{@Gil: pas tout \`a fait fin\dots.}
\[
%\hskip 2cm 
\E\, e^{\la V(X^{x_0}_{t})} \le e^{-\textcolor{myred}{\bar \a} t+ \la V(x_0)}+   \frac{\textcolor{myred}{\bar K}}{\textcolor{myred}{\bar \a}}\le  e^{\la V(x_0)}+   \frac{\textcolor{myred}{\bar K}}{\textcolor{myred}{\bar \a}}.
%\hskip 2cm \Box
\]
Now, one has for every $\Lambda\ge 0$,
 \begin{align*}
\int \big(e^{\la V(y)}\!\wedge \Lambda\big)  \hat \mu(\ud y) =\int \hat \mu(\ud x_0)\, \E ( e^{\la V(X_t^{x_0})}\wedge \Lambda)&\le \int\hat \mu(\ud x_0) \esp{ e^{\la V(X^{x_0}_t)}}\wedge \Lambda
\\
&\le  \int \hat \mu(\ud x_0) \Big(e^{-\textcolor{myred}{\bar \a} t+ \la V(x_0)}+   \frac{\textcolor{myred}{\bar K}}{\textcolor{myred}{\bar \a}} \big)\wedge \Lambda.
 \end{align*}
Letting successively $t\to +\infty$ and the $\Lambda\to +\infty$, first calling upon Lebesgue's dominated convergence theorem and then Fatou's Lemma and 
%Beppo Levi's 
monotone convergence theorem implies that
 \[
\hskip 5cm \int \big(e^{\la V(y)}\!\wedge \Lambda\big)  \hat \mu(\ud y)\le  \frac{\textcolor{myred}{\bar K}}{\textcolor{myred}{\bar \a}}.\hskip 5cm \Box
 \]

 \color{black}

 \subsection{Proof of Lemma \ref{le hoeffding like}} \label{subsec hoeffding like}

 \proof 
 $(a)$ We first observe that for any $z > 0$,
 \begin{align}
     \P\left(\int_0^t Z_s \ud W_s > \ell \psi(t)\right)
     =
     \P\left(e^{z\int_0^t Z_s \ud W_s} > e^{z\psi(t)\ell}\right).
 \end{align}
 Let $p$, $q>1$ be two H\" older conjugate exponents. Using exponential Markov inequality yields
%  \begin{align*}
%      \P\left(\int_0^t Z_s \ud W_s > \ell \psi(t)\right)
%      &\le \esp{e^{z\int_0^t Z_s \ud W_s}}e^{-z \psi(t)\ell}
%      \\
%      &=  \esp{e^{z\int_0^t Z_s \ud W_s-\tfrac p2\int_0^t |Z_s|^s\ud s} e^{\tfrac p2\int_0^t |Z_s|^s\ud s} }e^{-z \psi(t)\ell}
%      \\
%      &\le  \bigg(\esp{e^{pz\int_0^t Z_s \ud W_s-\tfrac{p^2}2\int_0^t |Z_s|^s\ud s}}\bigg)^{1/p} \bigg(\esp{e^{\tfrac{pq}{2}\int_0^t |Z_s|^2\ud s} }\bigg)^{1/q}e^{-z\psi(t)\ell}
%  \\
%      & \le \esp{e^{\frac{z^2}{2}p t|Z|^2_{\infty}}}e^{-z\psi(t)\ell}
%         \end{align*}
\color{myred}
\begin{align*}
     \P\left(\int_0^t Z_s \ud W_s > \ell \psi(t)\right)
     &\le \esp{e^{z\int_0^t Z_s \ud W_s}}e^{-z \psi(t)\ell}
     \\
     &=  \esp{e^{z\int_0^t Z_s \ud W_s-\frac{z^2}2\int_0^t |Z_s|^2\ud s} e^{\frac{z^2}2\int_0^t |Z_s|^2\ud s} }e^{-z \psi(t)\ell}
          \\ 
     &\le  \bigg(\esp{e^{z\int_0^t Z_s \ud W_s-\tfrac{z^2}2\int_0^t |Z_s|^2\ud s}}\bigg)   e^{\tfrac{{z^2}}{2} t |Z|_\infty^2}  e^{-z\psi(t)\ell}
 \\
     & \le e^{\tfrac{{z^2}}{2} t |Z|_\infty^2}  e^{-z\psi(t)\ell} %\color{black}
%      \\
%      &\le  \bigg(\esp{e^{pz\int_0^t Z_s \ud W_s-\tfrac{pz^2}2\int_0^t |Z_s|^2\ud s}}\bigg)^{1/p} \bigg(\esp{e^{\tfrac{q{z^2}}{2}\int_0^t |Z_s|^2\ud s} }\bigg)^{1/q}e^{-z\psi(t)\ell}
%  \\
%      & \le \esp{e^{\frac{z^2}{2} qt|Z|^2_{\infty}}}^{\frac1q}e^{-z\psi(t)\ell}
        \end{align*} \color{black}
        where the  expectation in the third line is $1$ owing to Novikov criterion.
   We obtain, by letting $p\to 1$ 
\begin{equation}
\P\left(\int_0^t Z_s \ud W_s > \ell \psi(t)\right) \le e^{\frac{z^2}2 t|Z|^2_\infty-z\psi(t)\ell}.
\end{equation}
%Optimizing the upper bound on $z$ , yields 
For $z=  \frac{t^{-1}\psi(t)\ell}{|Z|_\infty^2}$, one gets
 \begin{align}
     \P\left(\int_0^t Z_s \ud W_s > \ell \psi(t)\right)
     \le e^{-\frac{\ell^2}{2|Z|_{\infty}^2}t^{ - 1}\psi(t)^2 }.
 \end{align}
 Similar arguments leads to 
 \begin{align}
     \P\left(\int_0^tZ_s \ud W_s < -\ell  \psi(t)\right) \le 
     e^{-\frac{\ell^2}{2|Z|_{\infty}^2}t^{ - 1}\psi(t)^2 },
 \end{align}
 which concludes the proof for the first statement.
 
 \medskip
 %\color{blue}
\noindent  $(b)$ Using Markov inequality, we compute 
%\textcolor{red}{@Gill \1 JFC: why $a+2$? why not $a$ ?}
 \begin{align*}
     \P\left(\int_0^tZ_s \ud W_s >  \ell \psi(t)\right)
     &\le \frac{\esp{|\int_0^tZ_s \ud W_s|^{a+2}}}{\ell^{a+2}\psi(t)^{a+2}}
     %\\ &
     \le C^{\rm BDG}_{a+2} \frac{\esp{(\int_0^tZ_s^2 \ud s)^{\frac{a+2}2}}}{\ell^{a+2}\psi(t)^{a+2}}
 \end{align*}
 where we used Burkholder-Davis-Gundy  inequality for the last upper bound. Now, applying Jensen inequality, we get 
 \begin{align*}
     \P\left(\int_0^tZ_s \ud W_s >  \ell \psi(t)\right)
     &\le C \frac{t^{{\frac{a+2}2}}\esp{\frac1t \int_0^tZ_s^{a+2} \ud s}}{\ell^{a+2}\psi(t)^{a+2}}
     %\\&
     \le C \frac{t^{{\frac{a+2}2}}}{\ell^{a+2}\psi(t)^{a+2}},
 \end{align*}
 where we used the $L^{2+a}$ integrability assumption on $Z$ to get the last inequality.
 \\
 Using similar arguments, one proves
 \begin{align*}
     \P\left(\int_0^tZ_s \ud W_s < - \ell \psi(t)\right)
     %&\le C \frac{t^{{\frac{a+2}2}}\esp{(\int_0^tZ_s^{a+2} \ud s)}}{\ell^{a+2}\psi(t)^{a+2}}
    % \\&
    \le C \frac{t^{{\frac{a+2}2}}}{\ell^{a+2}\psi(t)^{a+2}},
 \end{align*}
 which concludes the proof for the second statement.
 \eproof

 \subsection{Various C\'esaro and Kronecker Lemmas}

 %\begin{Lemma}\label{le cesaro}
 %    Let $\rho \in [0,1)$ and $a>0$, then for $s \mapsto U(s)$ continuous and such that $\lim_{ s \rightarrow \infty} U(s) = \ell$,
 %    \begin{align}
 %        \frac{1-\rho}{t^{1-\rho}}\int_a^t s^{-\rho}U(s) \ud s \underset{t \rightarrow +\infty }{\longrightarrow }\ell \;.
 %    \end{align}
 %    Let $\alpha>0$, then 
 %    \begin{align}
 %        \alpha e^{-\alpha t} \int_0^t e^{\alpha u} U(s) \ud s \underset{t \rightarrow +\infty }{\longrightarrow } \ell.
 %    \end{align}
 %\end{Lemma}
 
   {
     \color{black}
     \begin{Lemma}[Continuous time C\'esaro's Lemma]\label{le Cesaro gil}
     Let $g:\R_+\to (0,+\infty)$ in $L^1_{loc}(\R_+,d\lambda)$ such that  $\int_0^{+\infty}g(t)dt =+\infty$. Then for every bounded Borel function  $U:\R_+\to \R $ such that $\lim_{t\to+\infty} U(t)= \ell\!\in \R$,  one has
     \begin{align}\label{eq  le Cesaro gil}
         \frac{\int_0^tg(s)U(s)ds}{\int_0^tg(s)ds}\ \underset{t \rightarrow +\infty }{\longrightarrow }\ell.
     \end{align}
 \end{Lemma}
 }
 {
     \color{black}
     \proof We may assume w.l.o.g. that $\ell=0$. Let $\ve>0$ and $A_{\ve}$ be such that $|U(t)|\le \ve $ for every $t\ge A_{\ve}$. Then 
     \[
     \bigg| \frac{\int_0^tg(s)U(s)ds}{\int_0^tg(s)ds}\bigg|\le \|U\|_{\sup}\frac{\int_0^{A_{\ve}}g(s)ds}{\int_0^tg(s)ds}+ \ve 
     \] 
 so that 
 $$
 \varlimsup_{t\to +\infty}\bigg| \frac{\int_0^tg(s)U(s)ds}{\int_0^tg(s)ds}\bigg|\le \ve\quad\mbox{for for every} \quad \ve >0.
 $$ 
    \eproof
 }
 
 %
 %\begin{Lemma}\label{le kronecker}
 %    Let $\rho \in [0,1)$ and $a>0$, then for $s \mapsto u(s)$ continuous and such that $\int_a^{+\infty} \frac{u(s)}{s} \ud s<+\infty$, 
 %    \begin{align}\label{eq re le kronecker}
 %        \frac{1}{t^{1-\rho}}\int_a^t s^{-\rho}u(s) \ud s \underset{t \rightarrow +\infty }{\longrightarrow }0 \;.
 %    \end{align}
 %\end{Lemma}
 %
 %{
 %    \color{blue}
 %    \proof 
 %        Let $U(s) = \int_s^{+\infty} \frac{u(s)}{s}\ud s$ for $s\ge a$. We observe that 
 %        \begin{align}
 %            t^{1-\rho}U(t) = U(1) + \int_1^t s^{1-\rho} \dot{U}(s) \ud s + \int_1^t s^{-\rho}{U}(s) \ud s,
 %        \end{align}
 %        leading to 
 %        \begin{align}
 %            \frac1{t^{1-\rho}}\int_1^t s^{-\rho} u(s) \ud s = \frac1{t^{1-\rho}}U(1)
 %            + \frac1{t^{1-\rho}}\int_1^t s^{-\rho}{U}(s) \ud s. 
 %        \end{align}
 %        Since $U(s)\rightarrow 0$, we deduce \eqref{eq re le kronecker} using Lemma \ref{le cesaro}.
 %    \eproof
 %}

 \begin{Lemma}[Kronecker's Lemmas]\label{le kronecker gil}
 Let $ a\ge 0$ and let $g:[a,+\infty)\to (0,+\infty)$ be a continuous  function such that  $G(t)= \int_a^{t}g(s)\ud s $ is well-defined for all $t \ge a$ and $ G(t) \to +\infty$ as $t\to +\infty$.

 \smallskip 
 \noindent $(a)$ {\em Continuous time regular Kronecker Lemma}.  Let $u:[a,+\infty]\to \R$ be a  continuous function. 
     Assume 
     $$
  \int_a^{t} (\log G)'(s) u(s) \ud s\underset{t \rightarrow +\infty }{\longrightarrow } I_\infty. %\; \R\mbox{-valued random variable}.
     $$
     Then 
     \begin{align}\label{eq re le kronecker gil}
         \frac{\int_a^t g(s)u(s)\ud s }{\int_a^t g(s)\ud s}\underset{t \rightarrow +\infty }{\longrightarrow }0 \;.
     \end{align}
     In particular if $\rho\!\in [0,1)$ and $\int_0^{t} \frac{u(s)}{s}\ud s\to \ell \in \R$ as $t\to +\infty$, then 
     \begin{align}\label{eq re le kronecker}
         \frac{1}{t^{1-\rho}}\int_0^t s^{-\rho}u(s) \ud s \underset{t \rightarrow +\infty }{\longrightarrow }0 \;.
     \end{align}
     
     \noindent $(b)$ {\em Stochastic Kronecker Lemma}. Let $W$ be an $(\cF_t)_t$-adapted standard Brownian motion  and $(u_t)_{t\ge a}$ be  an $(\cF_t)_{t\ge a}$-adapted process such that $ \int_a^{t} \big((\log G)'(s) u(s)\big)^2 \ud s<+\infty$ $\P$-$a.s.$. Assume 
     $$
  \int_a^{t} (\log G)'(s) u(s) \ud W_s \to I_\infty, \,\R  \mbox{-valued finite random variable as} \quad t\to +\infty.
     $$ 
 Then 
     \begin{align}\label{eq re le kronecker stoch}
         \frac{\int_a^t g(s)u(s)\ud W_s }{\int_a^t g(s)\ud s}\underset{t \rightarrow +\infty }{\longrightarrow }0 \;.
     \end{align}
     In particular if $\rho\!\in [0,1)$ and $\int_0^{t} \frac{u(s)}{s}\ud W_s\to I_\infty$, $\R$-valued finite random variable as $t\to +\infty$, then 
     \begin{align}\label{eq re le kronecker stoch}
         \frac{1}{t^{1-\rho}}\int_0^t s^{-\rho}u(s) \ud W_s \underset{t \rightarrow +\infty }{\longrightarrow }0 \;.
     \end{align}
 \end{Lemma}
 
 {
     \color{black}
     \proof 
         $(a)$ Let $U(t) = \int_a^t (\log G)'(s) u(s) \ud s $ for $t\ge a$. We observe that 
         \begin{align*}
            G(t)U(t) &= G(a)U(a) + \int_a^t U(s)g(s) \ud s + \int_a^t G(s)U'(s) \ud s,\\
            &=  \int_a^t U(s)g(s) \ud s + \int_a^t g(s) u(s)\ud s.
         \end{align*}
         Dividing the second equality by $G(t)$ and using the above C\'esaro Lemma~\ref{le Cesaro gil} yields
         \[  
           \frac{\int_a^t g(s)u(s)\ud s }{\int_a^t g(s)\ud s} = U(t)  - \frac{\int_a^tU(s)g(s)\ud s}{\int_a^t g(s)\ud s} \underset{t \rightarrow +\infty }{\longrightarrow } I_\infty-I_\infty = 0.
          \]
          
          \noindent The second claim follows by considering $g(t) = t^{-\rho}$. 
          \\
  $(b)$ Setting $U(t)=\int_a^t (\log G)'(s) u(s) \ud W_s$, the proof is the same as the one for $(a)$ since the integration by part formula is formally the same      owing to the locally  finite variation property of the function $G$.  
 %         
 %        leading to 
 %        \begin{align}
 %            \frac1{t^{1-\rho}}\int_1^t s^{-\rho} u(s) \ud s = \frac1{t^{1-\rho}}U(1)
 %            + \frac1{t^{1-\rho}}\int_1^t s^{-\rho}{U}(s) \ud s. 
 %        \end{align}
 %        Since $U(s)\rightarrow 0$, we deduce \eqref{eq re le kronecker} using Lemma \ref{le cesaro}.
     \eproof
 }